\documentclass[12pt]{amsart}
\usepackage{amsfonts}
\usepackage{graphicx}
\usepackage{epstopdf}
\usepackage{amsmath}
\usepackage{amsthm}
\usepackage{latexsym,bm}
\usepackage{amssymb}
\usepackage{esint}
\usepackage{enumitem}
\usepackage{indentfirst}
\usepackage{amscd}
\newtheorem{thm}{Theorem}[section]
\newtheorem{cor}[thm]{Corollary}
\newtheorem{lem}[thm]{Lemma}

\newtheorem{defn}[thm]{Definition}

\newtheorem{exm}{Example}
\newtheorem{hyp}{Hypothesis}

\newtheorem{Remark}{Remark}
\numberwithin{equation}{section}
\numberwithin{Remark}{section}

\begin{document}

\title{Constant $Q$-curvature metrics near the Hyperbolic metric}

\author{Gang Li}
\address{Department of Mathematics\\
       University of Notre Dame\\
        295 Hurley Hall \\
        Notre Dame, IN 46556, USA}
\address{Department of Mathematics\\
         Nanjing University\\
         Nanjing 210093, China}

\email{gli3@nd.edu}

\maketitle

\begin{abstract}
Let $(M,\,g)$ be a Poincar$\acute{\text{e}}$-Einstein manifold with                              
a smooth defining function. In this note, we prove that there are
infinitely many asymptotically hyperbolic metrics with constant
$Q$-curvature in the conformal class of an asymptotically hyperbolic
metric close enough to $g$. These metrics
are parametrized by the elements in the kernel of the 
linearized operator of the prescribed constant $Q$-curvature   
equation. A similar analysis is applied to a class of fourth order
equations arising in spectral theory.                     
\end{abstract}

\section{Introduction}


In this note we will discuss the prescribed constant $Q$-curvature
problem for asymptotically hyperbolic manifolds. We obtain the  
existence of a family of constant $Q$-curvature metrics in a small
neighborhood of any Poincar$\acute{\text{e}}$-Einstein metric,
parametrized by elements in the null space of the linearized
operator $L$ in $(\ref{linearoperator})$. Much of the analysis   
follows from Mazzeo's microlocal analysis method for elliptic edge
operators. Results in this setting have been proved for the scalar
curvature equation, see \cite{ACF}.                                                      



For $n\,\geq\,4$, a natural conformal invariant and the                                  
corresponding conformal covariant operator are the $Q$-curvature and
the fourth order Paneitz operator. Let $\text{Ric}_g$ and
$\text{R}_g$ be the Ricci curvature and the scalar curvature of
$(M,\, g)$. The $Q-Curvature$ and the $Paneitz\,operator$ are
defined as follows,
\begin{align*}
Q_g\,=\,\left\{\begin{array}{l@{\quad\quad}l} -\frac{1}{12}(\Delta_gR_g\,-\,R_g^2\,\,+\,3|\text{Ric}_g|^2), & n\,=\, 4,\\ \\
-\,\frac{2}{(n-2)^2}|\text{Ric}_g|^2\,+\,\frac{n^3-4n^2+16n-16}{8(n-1)^2(n-2)^2}R_g^2\,-\,\frac{1}{2(n-1)}\Delta_gR_g,
& n\,\geq\,5.\end{array} \right.\\
P_g(\varphi)\,=\,\left\{\begin{array}{l@{\quad\quad}l} \Delta_g^2\varphi\,-\,\text{div}(\frac{2}{3}R_g\,g\,-\,2\,\text{Ric}_g)d\varphi, & n\,=\, 4,\\ \\
\Delta_g^2\varphi\,-\,\text{div}_g(a_nR_g\,g\,-\,b_n\,\text{Ric}_g)\nabla_g
\varphi\,+\,\frac{n-4}{2}Q_g \varphi,& n\,\geq\,5,\end{array}
\right.
\end{align*}
\noindent where $a_n\,=\,\frac{(n-2)^2+4}{2(n-1)(n-2)}$,
$b_n\,=\,\frac{4}{n-2}$, $\text{div}_gX\,=\,\nabla_iX^i$ for any
smooth vector field $X$, and $\varphi$ is any smooth function on
$M$.

Let $\tilde{g}\,=\,\rho g$, with $\rho$ a positive function on $M$,
so that
\begin{align*}
\rho\,=\,\left\{\begin{array}{l@{\quad\quad}l} e^{2u}, & n\,=\,
4,\\
u^{\frac{4}{n-4}}, & n\,\geq\,5.\end{array} \right.
\end{align*}
The $Q$-curvature has the following transformation,                        
\begin{align*}
P_gu\,+\,2Q_g\,=\,2Q_{\tilde{g}}e^{4u},\,\,n\,=\,4,\\
P_gu\,=\,\frac{n-4}{2}Q_{\tilde{g}}\,u^{\frac{n+4}{n-4}},\,\,n>4.
\end{align*}
Note that Paneitz operator satisfies the following conformal covariance   
property for $\varphi\,\in\,C^{\infty}(M)$,
\begin{align*}
P_{\tilde{g}}\,\varphi\,=\,e^{-4u}P_g\,\varphi,\,\,\,n\,=\,4,\\
P_{\tilde{g}}(\varphi)\,=\,u^{-\frac{n+4}{n-4}}P_g(u\,\varphi),\,\,\,n>4.
\end{align*}

We want to find a function $u$ so that the metric $\tilde{g}$
satisfies $Q_{\tilde{g}}\,=\,f$  for a given function $f$. For the prescribed $Q$-curvature              
problem on closed manifold $M$ of dimension four there are many
results, see \cite{Chang-Yang}, \cite{DjadliaMalchiodi},                                                 
\cite{Matt.Gursky}, \cite{M.Gursky}. In \cite{Ndiay} a boundary
value problem for this problem is solved. A flow approach is
performed in \cite{S. Brendle}, see also \cite{Chen and Xu}. For                                         
$n\,\geq\,5$, see \cite{DHL},
\cite{Ndiay1} and \cite{Xu-Yang}. 

There are some interesting results for complete non-compact
manifolds. For Euclidean space $\mathbb{R}^{n}$,                                                
$n\,\geq\,4$, see \cite{Lin} and \cite{Wei and Ye}. 
In \cite{GRA}, using shooting method, the authors proved that there are infinitely many complete metrics with constant $Q$-curvature in the conformal class of the Poincar$\acute{\text{e}}$ disk with dimension $n \geq 5$, which are radially symmetric ODE solutions to the initial value problem parametrized by distinct given initial data at the origin.   
It is not difficult to prove that similar results hold for $n = 4$. 
Mazzeo pointed out that there should be a more general result of this type. In this paper, we solve a perturbation problem in the setting of asymptotically hyperbolic metrics close to a Poincar$\acute{\text{e}}$-Einstein metric. 
To give a precise statement we first need some definitions.                                                                   
\begin{defn}
Let $M$ be a smooth manifold of dimensional $n$, with smooth
boundary $\partial M$ of dimension $n-1$. Let $g$ be a complete
metric on $M\,=\,Int(\overline{M})$. We say that $g$ is asymptotically hyperbolic if    
there exists a smooth function $x$ on $\overline{M}$, with the
property that $x\,>\,0$ in $M$, and $x=0$ on $\partial M$, so that
the metric $h\,=\,x^2g$ is well defined and smooth on
$\overline{M}$, and $|dx|_h\,\big|_{\partial M}=\,1$. Here $x$ is
called a defining function of $g$. Moreover, if
$h\,\in\,C^{k,\,\alpha}$, for some positive integer $k$, we say that
$g$ is asymptotically hyperbolic of order $C^{k,\,\alpha}$. If $g$
is also Einstein, we call $g$ a Poincar$\acute{\text{e}}$-Einstein
metric, and $(M,\,g)$ a Poincar$\acute{\text{e}}$-Einstein manifold.
\end{defn}

Let $(M^n,\, g)$ be an asymptotically hyperbolic manifold of
dimension $n$, with $x$ as its smooth defining function. Actually,
we can choose $x$ so that $|d\,x|_h\,=\,1$ in a neighborhood of
$\partial M$, see \cite{Graham}, and here for simple notation we
always choose a
defining function in this sense except in Section $4$. We will mainly focus on the   
asymptotic behavior of the metric near $\partial M$, which is a
local discussion. Let $y$ be local coordinates on $\partial M$. In a
neighborhood of $\partial M$ in $\overline{M}$, we introduce the
local coordinates in the following way:
$(x,\,y)\,\in\,[0,\,\varepsilon)\,\times\,\partial M$ represent the
point moving from the point on $\partial M$ with local coordinate                   
$y$, along the geodesic which is the integral curve of $\nabla_h x$                   
for a length $x$ in the metric $h$. In local coordinates $(x,\,y)$,
\begin{align*}
h\,=\,x^2g\,=\,dx^2\,+\,\displaystyle\sum_{i,j=1}^{n-1}h_{ij}dy^idy^j.                               
\end{align*}

For convenience, let $\tilde{g}\,=\,\rho g$, with $\rho$ a positive
function on $M$, so that
\begin{align*}
\rho\,=\,\left\{\begin{array}{l@{\quad\quad}l} e^{2u}, & n\,=\,
4,\\
(\,1\,+\,u)^{\frac{4}{n-4}}, & n\,\geq\,5.\end{array} \right.
\end{align*}
Let the operator $\mathcal {E}$ be defined by
\begin{align}
\mathcal
{E}(u)\,=\,\left\{\begin{array}{l@{\quad\quad}l}\,P_gu\,+\,2Q_g\,-\,2Q_{\tilde{g}}e^{4u},
&\text{for}\,\,n\,=\,4,\\
\,P_g(1\,+\,u)\,-\,\frac{n-4}{2}Q_{\tilde{g}}(1\,+\,u)^{\frac{n+4}{n-4}},
& \text{for}\,\,n\,\geq\,5.\end{array} \right.
\end{align}
To solve the prescribed $Q$-curvature problem amounts to finding a
solution to
\begin{align}\label{NLPS}
\mathcal {E}(u)\,=\,0.
\end{align}
We define the linear operator $L\,=\,L_g$ as follows,                                      

\begin{align}\label{linearoperator}                              
L(u)\,=\,\left\{\begin{array}{l@{\quad\quad}l} P_gu\,-\,8Q_{g}u, & n\,=\, 4,\\ \\
P_gu\,-\,\frac{n+4}{2}Q_{g}u, & n\,\geq\,5.\end{array} \right.
\end{align}

Let $(x,y)$ be the local coordinates of $M$ near the boundary
defined as above. Let $\mathcal {V}_e$ be the collection of the
smooth vector fields on $\overline{M}$, which restricted in the
neighborhood of $\partial M$, are generated by
$\{x\partial_x,\,x\partial_{y^1},\,...,x\partial_{y^{n-1}}\}$ with                            
smooth coefficients on $\overline{M}$.                                                        

Next we introduce the weighted spaces that we will be using. First,
the weighted Sobolev spaces,
\begin{align*}
x^{\delta}H_e^{m}(M,\,\Omega^{\frac{1}{2}})=\{\,u\,=\,x^{\delta}\,v:\,\,V_1...V_jv\in\,L^2(M,\,\Omega^{\frac{1}{2}}),\,\forall\,j\,\leq\,m,\,V_i\,\in\,\mathcal                           
{V}_e\,\},
\end{align*}
where $m\,\in\,\mathbb{N}$, $\delta\,\in\,\mathbb{R}$, and $\Omega^{\frac{1}{2}}\,=\sqrt{dx dy}$ is the half-density.                
We also introduce the weighted H$\ddot{\text{o}}$lder space, 
\begin{align*}
x^{\delta}\Lambda^{m,\,\alpha}=\,x^{\delta}\Lambda^{m,\,\alpha}(M,\,\Omega^{\frac{1}{2}})\,=\,\{\,u\,=\,x^{\delta}v\sqrt{dx\,dy}:\,\,V_1...V_jv\,\in\,\Lambda^{0,\,\alpha},\,\forall\, j\,\leq\,m,\,V_i\,\in\,\mathcal     
{V}_e\,\},
\end{align*}
with $m\,\in\,\mathbb{N}$, $\delta \,\in\,\mathbb{R}$, and $0 < \alpha <1$, where                           
$\Lambda^{0,\alpha}(M)$ is the space of half-densities
$u\,=\,v\sqrt{dx\,dy}$ such that
\begin{align*}
\|v\|_{\Lambda^{0,\alpha}(M)}\,=\,\text{sup}\,|v|\,+\,\text{sup}\frac{(x\,+\,\tilde{x})^{\alpha}\,|v(x,y)\,-\,v(\tilde{x},\tilde{y})|}{|x\,-\,\tilde{x}|^{\alpha}\,+\,|y\,-\,\tilde{y}|^{\alpha}}\,<\,\infty.   
\end{align*}
 We will use
the norm
\begin{align*}
\|u\|_{x^{\delta}\Lambda^{k,\,\alpha}(M)}\,=\,\displaystyle\sum_{m=0}^k\sum_{|\gamma|\,=\,m}\|\partial_{e}^{\gamma}v\|_{0,\,\alpha},   
\end{align*}
with $\partial_{e}\,\in\,\mathcal {V}_e$ and $u\,=\,x^{\delta}v$.                                  

In this paper, we always assume $n\geq 4$ to be the dimension of $M$.                         
With these definitions, we can now state our main result:
\begin{thm}\label{Hyperbolicmetrico}
Let $(B^n_1(0), g)$ be the Poincar$\acute{\text{e}}$ disk, of
dimension $n\geq 4$. Also, let $x$ be a smooth defining function of                       
$g$. Let $L$ be the linear operator defined in
$(\ref{linearoperator})$. Let $\nu$ be a constant in the interval                         
$(0,\,\frac{n-1}{2})$. Then,
\begin{enumerate}[label=\roman{*})]

\item\label{Theorem1_1} Kernel of $L$ in the weighted space                                                  
$x^{\nu}\Lambda^{4,\,\alpha}(M)$ for $0<\alpha<1$ is of infinite dimension. Also, $L$                     
is surjective.                                                                             
For each element $v$ in the kernel $\ker(L)$ with sufficiently small norm,                       
and a given function                                                                       
$Q_{\tilde{g}}\,\in\,\Lambda^{0,\,\alpha}(M,\,\sqrt{dx dy})$ so that
$(Q_{\tilde{g}}\,-\,Q_g)$ is in
$x^{\nu}\Lambda^{0,\,\alpha}(M,\,\sqrt{dx dy})$ with the norm $
\|Q_{\tilde{g}}\,-\,Q_g\|_{x^{\nu}\Lambda^{0,\,\alpha}}$ small
enough, there exists a unique solution
$u\,\in\,x^{\nu}\Lambda^{4,\,\alpha}(M)$ to the problem                                    
$(\ref{NLPS})$, so that the projection $P_1$ ( see in Theorem                               
\ref{HLSemiFredholm}) of $u$ onto $\ker(L)$ is given by $v$.

\item Moreover, if $Q_{\tilde{g}}\,=\,Q_g$, $u$ has the expansion near the boundary                         
\begin{align}\label{mainterms}
u(x,\,y)\,\sim\,(u_{00}(y)x^{\frac{n-1}{2}\,+
\,i\beta}\,+\,u_{10}(y)x^{\frac{n-1}{2}\,-\,i\beta})\,
+ \,o(x^{\frac{n-1}{2}}),                                                                                                       
\end{align}
with $\beta\,=\,\frac{\sqrt{n^2+2n-9}}{2}$ and $i\,=\,\sqrt{-1}$, where $u_{00}$ and $u_{10}$ are generally distributions of negative                    
order. Also, $u$ will have the following expansion with smooth coefficients,                           
\begin{align}\label{fomalexpansion}
u(x,\,y)\,\sim\,\displaystyle\sum_{j\,=\,0}^{+\,\infty}(u_{0j}(y)x^{\frac{n-1}{2}\,+\,i\beta\,+\,j}\,+\,u_{1j}(y)x^{\frac{n-1}{2}\,-\,i\beta\,+\,j}\,+\,u_{2j}(y)x^{n\,+\,j}),       
\end{align}
in the sense that
\begin{align*}
u(x,y)\,-\,\displaystyle\sum_{j=0}^k(u_{0j}(y)x^{\frac{n-1}{2}\,+\,i\beta\,+\,j}\,+\,u_{1j}(y)x^{\frac{n-1}{2}\,-\,i\beta\,+\,j})\,=\,o(x^{\frac{n-1}{2}+k}),    
\end{align*}
with $\beta\,=\,\frac{\sqrt{n^2+2n-9}}{2}$ for each $k\,\geq\,0$,                                                           
if $v\,=\,P_1u$ has an expansion of this form with smooth                                                                   
coefficients and $1\,\leq\,\nu\,<\,\frac{n-1}{2}$.
\end{enumerate}
\end{thm}
For kernel elements having an expansion with smooth coefficients,                                                            
one can prescribe the leading terms for them, see Remark
\ref{Remark_smoothness}.
\begin{Remark}
The ODE result in \cite{GRA} only gives existence of radially     
symmetric constant $Q$-curvature metrics in the conformal class of
the hyperbolic metric, but allows the metric to be far away from the hyperbolic metric. As a perturbation result, our theorem gives                               
the existence of solutions in the conformal class of metrics in a
small neighborhood of the
hyperbolic metric, more precisely, see Theorem \ref{generating}.                                  
\end{Remark}

Using boundary regularity results and the unique continuation 
property on the boundary, as a slight extension of the above theorem we have the            
following result. Note that both boundary regularity results and the 
unique continuation property approach need $x$ and $h\,=\,x^{2}g$ to                
be smooth enough on $\overline{M}$.
\begin{thm}\label{PCPE1}
Let $(M^n,g)$, $n\geq 4$, be a Poincar$\acute{\text{e}}$-Einstein manifold with                                                  
the defining function $x$ and the metric $h\,=\,x^2g$ smooth up to
the boundary. Suppose also that
$L\,:\,x^{\nu}\Lambda^{4,\,\alpha}(M)\,\rightarrow\,x^{\nu}\Lambda^{0,\,\alpha}(M)$,
where $0\,<\,\nu\,<\frac{n-1}{2}$ and $0<\alpha<1$, is defined in                                                  
(\ref{linearoperator}). Then,
\begin{enumerate}[label=\roman{*})]

\item\label{Theorem2_1} Kernel of $L$ in the weighted space 
$x^{\nu}\Lambda^{4,\,\alpha}(M)$ is of infinite dimension. Also, $L$
is surjective. For each element $v$ in the kernel with its norm
small enough, and a given function                                                                       
$Q_{\tilde{g}}\,\in\,\Lambda^{0,\,\alpha}(M,\,\sqrt{dx dy})$ so that
$(Q_{\tilde{g}}\,-\,Q_g)$ is in
$x^{\nu}\Lambda^{0,\,\alpha}(M,\,\sqrt{dx dy})$ with the norm $
\|Q_{\tilde{g}}\,-\,Q_g\|_{x^{\nu}\Lambda^{0,\,\alpha}}$ small
enough, there exists a unique solution $u$ to the problem
$(\ref{NLPS})$, so that the projection $P_1$( see in Theorem \ref{HLSemiFredholm}) of $u$ onto $\ker(L)$ is                              
given by $v$.

\item Moreover, if $Q_{\tilde{g}}\,=\,Q_g$, then $u$ has the                                       
expansion result as in Theorem \ref{Hyperbolicmetrico}.
\end{enumerate}
\end{thm}

Since this is a perturbation result, we first discuss the linear
problem. Using Mazzeo's approach in \cite{Mazzeo1}, we obtain the
semi-Fredholm property for the linear operator
$(\ref{linearoperator})$:
\begin{thm}\label{SLSemiFredholm}
Let $(M^n,g)$ be an asymptotically hyperbolic manifold with
defining function $x$ and the metric $h\,=\,x^2\,g$ smooth up to                          
the boundary, then the linear operator
$L\,:\,x^{\delta}H_e^{4}(M)\,\rightarrow\,x^{\delta}L^{2}(M,\,\sqrt{dx\,dy})$
as in (\ref{linearoperator}), is essentially injective if                                 
$\delta\,>\,\frac{n}{2}$ and $\delta\,\neq\,n+\frac{1}{2}$, with                          
infinite dimensional cokernel, and $L$ is essentially surjective if
$\delta\,<\,\frac{n}{2}$  and $\delta\,\neq -\frac{1}{2}$, with                           
infinite dimensional kernel.( Here essentially injective means that
the null space of $L$ is at most finitely dimensional, and
essentially surjective means that $L$ has closed range and with at
most finitely dimensional cokernel.) Moreover, in both cases, $L$
has closed range, and admits a generalized inverse $G$ and
orthogonal projectors $P_1$ onto the nullspace and $P_2$ onto
orthogonal complement of the range of $L$ which are edge operators,
such that,
\begin{align*}
GL\,=\,I\,-\,P_1,\\
LG\,=\,I\,-\,P_2.
\end{align*}
\end{thm}
The corresponding theorem for the weighted H$\ddot{\text{o}}$lder
space is as follows.
\begin{thm}\label{HLSemiFredholm}
Let $(M^n,g)$ be an asymptotically hyperbolic manifold with
defining function $x$ and the metric $h\,=\,x^2\,g$ smooth up to                          
the boundary. Let
$0<\alpha<1$. The linear operator                                                                         
$L\,:\,x^{\nu}\Lambda^{4,\,\alpha}(M)\,\rightarrow\,x^{\nu}\Lambda^{0,\,\alpha}(M)$                       
as in (\ref{linearoperator}), is essentially injective if                                                 
$\nu\,>\,\frac{n-1}{2}$ and $\nu\,\neq\,n$, with infinite                                                 
dimensional cokernel; and
$L$ is essentially surjective if $\nu\,<\,\frac{n-1}{2}$ and $\nu\neq -1$, with                           
infinite dimensional kernel. Moreover, in both cases, $L$ has closed
range. Also, $x^{\nu}\Lambda^{4,\,\alpha}(M)$ has the topological
splitting of the following direct sum
$x^{\nu}\Lambda^{4,\,\alpha}(M)\,=\,P_1(x^{\nu}\Lambda^{4,\,\alpha}(M))\,\oplus\,(I\,-\,P_1)(x^{\nu}\Lambda^{4,\,\alpha}(M))$,
which are the projection to the null space of $L$ and its
topological complement for the second case. Similarly as the theorem
with weighted Sobolev spaces, there is a corresponding splitting of
$x^{\nu}\Lambda^{0,\,\alpha}(M)$ for $\nu\,>\,\frac{n-1}{2}$.                                             
\end{thm}

The paper is organized as follows. In Section $2$, we study the
linear elliptic edge operator $L$ defined in
(\ref{linearoperator}), 
and obtain the semi-Fredholm property of the linear operator $L$. In
Section $3$, we obtain that if the linear operator $L$ with respect
to the initial asymptotically hyperbolic metric $g$ is surjective in                  
a suitable weighted H$\ddot{\text{o}}$lder space, there are
infinitely many solutions to the prescribed $Q$-curvature problem
with $Q_{\tilde{g}}$ a small perturbation of $Q_g$, and the                           
solutions are
parametrized by the elements in the kernel of $L$. Then we give the proof of Theorem \ref{Hyperbolicmetrico} and Theorem \ref{PCPE1}.      
Using a special weighted H$\ddot{\text{o}}$lder space, in Section $4$, we prove a perturbation result                   
for the prescribed constant $Q$-curvature problem for a
Poincar$\acute{\text{e}}$-Einstein metric. In Section $5$, we give a                   
similar discussion to the prescribed $U$-curvature
equations.\vskip5pt

\section {Semi-Fredholm properties of the linearized
operator}\label{section_a}
\vskip10pt

In the following, we will discuss the local parametrix for $L$ and
the Fredholm property of $L$. A clear feature is that the elliptic operators $L$ under consideration here are degenerate near infinity. Here we review some of the material developed by Mazzeo and others  
in the theory of elliptic edge operators.

As in the introduction, let $(M^n,\, g)$ be an asymptotically
hyperbolic manifold of dimension $n$,  with
defining function $x$ and the metric $h\,=\,x^2\,g$ smooth up to                          
the boundary. Let $(x,y)$ be the local                                                               
coordinates of $M$ near the boundary, and $\mathcal {V}_e$ be
defined in the introduction. The one-forms dual to the vector fields
which are elements in $\mathcal {V}_e$ are smooth one forms in $M$,                     
restricted on the neighborhood of $\partial M$ generated linearly by
$\{\frac{d\,x}{x},\,\frac{d\,y^1}{x},\,...,\frac{d\,y^{n-1}}{x}\}$                                                              
with
coefficients smooth up to $\partial M$. Generally, a left or right parametrix $E$ of an elliptic operator $L$ on $M$                    
is a pseudo-differential operator with the property that                                                                                
\begin{align*}
EL\,=\,\text{Id}\,+\,R_1,\,\,\text{or}\,\,LE\,=\,\text{Id}\,+\,R_2,
\end{align*}
with $R_1,\,R_2$ compact operators.

The Schwartz kernel of an interior parametrix of the linear
operator $L$ is a distribution on $M \times M$, and for "interior" we mean that the parametrix has singularity near the boundary which will be 
explained in the following. Let $(x,\,y)$ and
$(\tilde{x},\,\tilde{y})$ be local coordinates on each copy of $M$
near the boundary. We know that the parametrix is smooth, except for
the singularity along the diagonal
$\Delta\,=\,\{x\,=\,\tilde{x},\,y\,=\,\tilde{y}\}$, as in the case
of compact manifolds. Moreover, here due to the degeneration of the
edge operator $L$, as $x,\,\tilde{x}\,\to\,0$, we also have the
important additional singularity at the intersection of $\Delta$ and
the corner, which is
$S\,=\,\{\,x=\,\tilde{x}\,=\,0,\,y\,=\,\tilde{y}\}$. To deal with
the boundary singularity, we introduce a new manifold
$M_0^2\,=\,M\,\times_0\,M$, by blowing-up $M\,\times\,M$ along $S$. 
Actually, if we use polar coordinates for $M\,\times\,M$ near the
corner,
\begin{align*}
&r\,=\,(x^2\,+\,|y\,-\,\tilde{y}|^2\,+\,\tilde{x}^2)^{1/2}\,\in\,\mathbb{R}^{+},\\
&\Theta\,=\,(x,\,y\,-\,\tilde{y},\,\tilde{x})/r\,\in\,S^n_{++}\,=\,\{\Theta\,\in\,S^n,\,\Theta_0,\,\Theta_n\,\geq\,0\},
\end{align*}
we know that the level set of $r\,=\,R$ is a submanifold of
dimensional $2n-1$ for $R\,>\,0$, while $S\,=\,\{\,r\,=\,0\}$ is
singular. More precisely, let $M_0^2$ be the lift of $M\,\times\,M$                                                                            
such that
it is the same as $M\,\times\,M$ away from $S$, but near the corner, 
it is represented by the lift of the polar coordinates, smoothly.
Hence, $S_{11}\,=\,\{\,r=0\}$ is a $(2n-1)$-dimensional submanifold                                                                           
of $M_0^2$. Let $b$ be the natural projection map from $M_0^2$ to
$M\,\times\,M$. For the convenience of calculation, as in
\cite{Mazzeo1}, we introduce two systems of local coordinates on   
$M_e^2$, $(s,v,\tilde{x},\tilde{y})$ and $(x,y,t,w)$,              
where
\begin{align*}
s\,=\,x/\tilde{x},\,v\,=\,\frac{y-\tilde{y}}{\tilde{x}};\,\,t\,=\,\tilde{x}/x,\,w\,=\,\frac{\tilde{y}-y}{x}.                 
\end{align*}
Changing variables in these two coordinates, 
\begin{align*}
x\partial_x\,=\,s\partial_s\,=\,x\partial_x\,-\,w\partial_w\,-\,t\partial_t,\,\,
\text{and}\,\,x\partial_y\,=\,s\partial_v\,=\,x\partial_y\,-\,\partial_w.                                                 
\end{align*}
In the following with out loss of generality we only need to
consider $(s,v,\tilde{x},\tilde{y})$. Viewing elements in                                                              
$\mathcal{V}_e$ as first order differential operators, we denote
$\text{Diff}_e^{*}(M)$ the algebra generated by $\mathcal {V}_e$
with coefficients in the ring $C^{\infty}(\overline{M})$, and with
the product given by composition of operators. Let
$\text{Diff}_e^m(M)$ be the linear subspace of differential             
operators which are of $m$-th order. Then for
$L\,\in\,\text{Diff}_e^m(M)$, it has the form
\begin{align}\label{edgeelliptic}                               
L\,=\,\sum_{j+|\alpha|\,\leq\,m}\,a_{j,\alpha}(x,y)(x\partial_x)^j(x\partial_y)^{\alpha},                               
\end{align}
with $a_{j,\alpha} \,\in\,C^{\infty}(\overline{M})$, in the coordinate chart $(x, y)$. The symbol of $L$ is                  
\begin{align*}
\sigma_e(L)(x,y;\xi,\eta)\,=\,\sum_{j+|\alpha|\,=\,m}\,a_{j,\,\alpha}(x,y)\xi^j\eta^{\alpha}.
\end{align*}
$L$ is elliptic if $\sigma_e(L)(x,y;\xi,\eta)\,\neq\,0$, for
$(\xi,\,\eta)\,\neq\,0$. It is easy to check that $\Delta_g$ and the
linear operator $L$ in (\ref{linearoperator}) are elliptic. $L$ in
$(\ref{edgeelliptic})$ can be considered as a lift to $M_e^2$ as follows,   
\begin{align*}
L\,=\,\displaystyle \sum_{j+|\alpha|\,\leq\,
m}\,a_{j,\alpha}(x,y)(x\partial_x)^j(x\partial_y)^{\alpha}\,=\,\displaystyle
\sum_{j+|\alpha|\,\leq\,\,
m}\,a_{j,\alpha}(s\tilde{x},\,\tilde{y}\,+\,\tilde{x}v)(s\partial_s)^j(s\partial_v)^{\alpha}.                
\end{align*}
Let $N(L)$ be the normal operator of $L$, so that                                                            
\begin{align*}
N(L)\,=\,\displaystyle \sum_{j+|\alpha|\leq
m}\,a_{j,\alpha}(0,\tilde{y})(s\partial_s)^j(s\partial_v)^{\alpha},                                       
\end{align*}
is the restriction to $S_{11}$ of the lift of $L$ to $M_e^2$. The                                      
normal operator is an important approximation of $L$ near the
boundary. 
For the linear operator $L$ in (\ref{edgeelliptic}), 
\begin{align*}
L\phi\,=\,\displaystyle \sum_{j+|\alpha|\leq
m}\,a_{j,\alpha}(0,\,y)(x\partial_x)^j(x\partial_y)^{\alpha}\phi\,+\,E\phi,                              
\end{align*}
any smooth function $\phi$, with the error term                                        
\begin{align*}
E\phi\,=\,x\,\displaystyle\sum_{j+|\alpha|\leq
m}b_{j,\alpha}(x,y)(x\partial_x)^j(x\partial_y)^{\alpha}\phi,                            
\end{align*}
for $x>0$ small, with the                               
coefficients $b_{j,\alpha}$ smooth up to the boundary.
\begin{defn}
The indicial family $I_{\zeta}(L)$ of $L \in \text{Diff}_e^k(M)$ is                    
defined to be the family of operators
\begin{align*}
L(x^{\zeta}(\log(x))^pf(x,y))\,=\,x^{\zeta}(\log(x))^pI_{\zeta}(L)f(0,y)\,+\,O(x^{\zeta}(\log(x))^{p-1}),
\end{align*}
for $\,f\,\in\,C^{\infty}(M),\,
\zeta\,\in\,\mathbb{C},\,p\,\in\,\mathbb{N}_0$.

\end{defn}
There exists a unique dilation-invariant operator $I(L)$, which is
called the indicial operator, such that
\begin{align*}
I(L)(y,\,s\partial_s)s^{\zeta}f(y)\,=\,s^{\zeta}I_{\zeta}(L)f(y).
\end{align*}
In local coordinates near the boundary,
$I(L)\,=\,\sum_{j\,\leq\,k}a_{j,\,0}(0,y)(s\partial_s)^j$.

\begin{defn}
If $L\,\in\,\text{Diff}_e^*(M)$ is elliptic, we denote $spec_b(L)$                                         
as the boundary spectrum of $L$, which is the set of
$\zeta\,\in\,\mathbb{C}$, for which $I_{\zeta}(L)\,=\,0$.
\end{defn}
Let $(M,\,g)$, $x$, and $h$ be defined as above. Denote $S_x$ as the
level set of $x$ ( also denoted as $x_0$ for convenience), and the
coordinates $(y_1,...,y_{n-1})\,=\,y$. We now use
this point of view to deal with our linearized operator $(\ref{linearoperator})$.                               

In a neighborhood of $\partial M$, we have the following,                                                                 
\begin{align}\label{Ricciformula}
\text{Ric}_g\,=\,\text{Ric}_h\,+\,x^{-1}[(n-2)\text{Hess}_hx\,+\,\Delta_hx\,h]\,-\,(n-1)x^{-2}|dx|_h^2h,
\end{align}
and
\begin{align}\label{scalarcurvatureform}
R_g\,=\,-n(n-1)|dx|_h^2\,+\,(2n-2)x(\Delta_hx)\,+\,x^2R_h,
\end{align}
where $|dx|_h\,=\,1$, and                      
\begin{align*}
(\text{Hess}_h)_{ij}(x)\,=\,\nabla_i^h\nabla_j^h(x)\,=\,\partial_i\,\partial_j(x)\,-\,\Gamma_{ij}^s\partial_s(x)\,=\,-\,\Gamma_{ij}^0\,=\,\frac{1}{2}\partial_x
h_{ij}\,=\,B_{ij},
\end{align*}
with $B_{ij}$ the second fundamental form of $S_x$,                               
for $i,\,j\,>0$; and $(\text{Hess}_h)_{ij}(x)\,=\,0$
otherwise. Also $\Delta_hx\,=\,\text{tr}_h(\text{Hess}_h)=\,H(h)$,                                                          
with $H(h)$ the mean curvature of the level set of $x$ in the metric
$h$. Here $\Gamma_{ij}^k$ is the Christoffel symbol with respect to
$h$. Note that $\Delta_g$ in our paper is the trace of                                
$\text{Hess}_g$, with negative eigenvalues:
\begin{align}\label{Laplacian}
\Delta_g\,u\,&=\,g^{ij}(\partial_i\partial_j\,-\,\Gamma_{ij}^k\partial_k)\,u\\
&=\,x^2\Delta_h\,u\,+\,(2-n)\,x\,(\nabla_h x,\, du)\\
&=\,(2-n)\,x\partial_x\,u\,+\,x^2(\partial_x^2u\,+\,\Delta_y\,u\,+\,H(h)\partial_x\,u),                                                                   
\end{align}
where $\Delta_y$ is the Laplacian on the level set $S_x$ of $x$, in
the
induced metric $h\big|_{S_x}$.                                                                                                           

Near the boundary, the $Q -$curvature is
\begin{align*}
Q_g\,&=\,-\,\frac{2}{(n-2)^2}(n-1)^2n\,+\,\frac{n^3-4n^2+16n-16}{8(n-1)^2(n-2)^2}n^2(n-1)^2\,+\,O(x)\\
&=\,\frac{n(n^2-4)}{8}\,+\,O(x),                                                                                                        
\end{align*}
for $n\,\geq\,5$, and $Q_g\,=\,3\,+\,O(x)$, for $n\,=\,4$.

In the following of this section we will discuss about the linear
operator $L$ in (\ref{linearoperator}).                                                     
Note that
\begin{align*}
L\,\phi\,&=\,\Delta_g^2\phi\,-\,\text{div}_g(a_nR_g\,g\,-\,b_n\,\text{Ric}_g)\nabla_g
\phi\,-\,4\,f\,\phi\\
&=\,\Delta_g^2\phi\,-\,a_nR_g\Delta_g\phi\,+\,b_n\,\text{Ric}_{ij}^g\nabla_g^i\nabla_g^j
\phi\,-a_n(\nabla_gR_g,\,\nabla_g\phi)\,+\,b_n\,\nabla^i_g\text{Ric}_{ij}\nabla_g^j\phi-\,4\,f\,\phi,\\
&=\,\Delta_g^2\phi\,-\,a_nR_g\Delta_g\phi\,+\,b_n\,\text{Ric}_{ij}^g\nabla_g^i\nabla_g^j
\phi\,+\,(-a_n\,+\,\frac{b_n}{2})(\nabla_gR_g,\,\nabla_g\phi)-\,4\,f\,\phi\\
&=\,\Delta_g^2\phi\,-\,a_nR_g\Delta_g\phi\,+\,b_n\,\text{Ric}_{ij}^g\nabla_g^i\nabla_g^j
\phi\,+\frac{6-n}{2(n-1)}(\nabla_gR_g,\,\nabla_g\phi)-\,4\,f\,\phi,
\end{align*}
with $f\,=\,Q_g$ for $n\,\geq\,5$, and $f\,=\,2\,Q_g$ for $n\,=\,4$.
For the third equality, we use the second Bianchi identity. Also, 
\begin{align*}
\Delta_g\phi\,=\,x^2\Delta_h\phi\,-\,(n-2)\,x\,(\nabla_hx,\,d\phi)_h\,=\,x^2\,\Delta_h\phi\,-\,(n-2)x\partial_x\phi,
\end{align*}
and
\begin{align*}
R_{ij}(g)\,\nabla_g^i\,\nabla_g^j\,\phi\,&\sim\,[-(n-1)x^2h_{ij}\,+\,O(x^3)]x^{-4}\nabla^i\,\nabla_g^j\,\phi\,\\             
&=\,-(n\,-\,1)(\Delta_g\,\phi\,+\,O(x)\,p(x,\,y,\,x\partial_x,\,x\partial_y)\phi),
\end{align*}
for some smooth function $p(\cdot)$. As a consequence,
\begin{align*}
L\,\phi\,=\,&\Delta_g^2\phi\,-\,a_nR_g\Delta_g\phi\,+\,b_n\,\text{Ric}_{ij}^g\nabla_g^i\nabla_g^j
\phi\,+\frac{6-n}{2(n-1)}(\nabla_gR_g,\,\nabla_g\phi)-\,4\,f\,\phi\\
=\,&\Delta_g^2\phi\,-\,a_n(-n(n-1)\,+\,O(x))\Delta_g\phi\,+\,b_n\,(-(n-1)\Delta_g\,\phi\\
&+\,O(x)p(x,\,y,\,x\partial_x,\,x\partial_y)\phi)
+\frac{6-n}{2(n-1)}(-(2n-2)x^2H(h|_{S_x})\partial_x\phi\\
&+\,O(x^3)|\nabla_y\phi|)-\,(\frac{1}{2}n(n^2\,-\,4)\,+\,O(x))\phi,
\end{align*}
and then, by definition,                           
\begin{align*}
N(L)\,=\,[(s\partial_s)^2\,-\,(n-1)s\partial_s\,+\,s^2\Delta_v\,-\,n][(s\partial_s)^2\,-\,(n-1)s\partial_s\,+\,s^2\Delta_v\,+\,\frac{n^2-4}{2}].      
\end{align*}
In addition,                                       
\begin{align*}
I(L)\,=\,\big((s\partial_s)^2\,-\,(n-1)s\partial_s\,-\,n\big)\big((s\partial_s)^2\,-\,(n-1)s\partial_s\,+\,\frac{n^2-4}{2}\big).
\end{align*}
Let $\phi\,=\,s^{\zeta}$, and $I(L)\phi\,=\,0$. Solving the
equation, we get the indicial roots $\zeta$, given by                                                                                    
\begin{align*}
spec_b(L)\,=\,\{n,\,-1,\,\frac{n-1}{2}\,-\,i\,\frac{\sqrt{n^2\,+\,2n\,-\,9}}{2},\,\frac{n-1}{2}\,+\,i\,\frac{\sqrt{n^2\,+\,2n\,-\,9}}{2}\}.
\end{align*}
\noindent Let $\Lambda$ be the indices set
\begin{align}\label{indicesset}
\Lambda\,=\,\{\,\frac{1}{2}\,+\,\text{Re}(\delta);\,\delta\,\in\,spec_b(L)\}.
\end{align}
The operator $N(L)$ acts on functions defined on
$\mathbb{R}_s^{+}\,\times\,\mathbb{R}_v^{n-1}$ for each fixed                                             
$\tilde{y}$, with coordinates $(s,\,v)$. For our linear operator                                          
$L$, $N(L)$ does not depend on $\tilde{y}$. We now take the Fourier
transformation of $N(L)$ in $v$ direction,                                                                
\begin{align*}
\widehat{N(L)}\,=\,\displaystyle
\sum_{j\,+\,|\alpha|\,\leq\,m}a_{i,\alpha}(s\partial_s)^j(i\,s\,\eta)^{\alpha}.
\end{align*}
We have the symmetry of dilation:
\begin{align*}
a_{j\,\alpha}\,(s\,\partial_s)^j\,(s\,\partial_y)^{\alpha}\,=\,a_{j\,\alpha}\,(ks\,\partial_{k\,s})^j\,(k
s\,\partial_{k\,y})^{\alpha},
\end{align*}
for any $k\,\in\,\mathbb{R}\,-\,\{0\}$. Let $t\,=\,s\,|\eta|$, then
\begin{align*}
\widehat{N(L)}(s,\,\eta)\,=\,\displaystyle
\sum_{j\,+\,|\alpha|\,\leq\,m}a_{i,\alpha}(0,\,\tilde{y})(t\partial_t)^j(i\,t\,\widehat{\eta})^{\alpha},
\end{align*}
 which is denoted as $L_0(t,\,\hat{\eta})$, where
$\hat{\eta}\,=\,\frac{\eta}{|\eta|}$. This is a family of totally
characteristic operators on $\mathbb{R}_{+}^n$ and generally its
coefficients depend on $\tilde{y}$. Now we have fixed
$\widehat{\eta}$ in the formula, and it has no scaling freedom in
this direction.

Let $\mathcal {H}^{m,\delta,l}$ be the weighted Sobolev space
\begin{align*}
\mathcal{H}^{m,\,\delta,\,l}\,=\,\{f:\,\phi(t)f\,\in\,t^{\delta}H_e^m(\mathbb{R}^{+}),\,(1\,-\,\phi(t))f\,\in\,t^{-l}H^m(\mathbb{R}^{+})\},
\end{align*}
with $\phi\,\in\,C_0^{\infty}(\mathbb{R}^{+})$, and $\phi(t)\,=\,1$
in a neighborhood of $t\,=\,0$. Note that                              
\begin{align*}
L_0:\,t^{\delta}\mathcal{H}^{m,\,\delta,\,l}\,\to\,t^{\delta}\mathcal
{H}^{m-4,\,\delta,\,l+4}
\end{align*}
is bounded.

For our linear operator $L$,
\begin{align*}
\widehat{N(L)}=[(s\partial_s)^2-(n-1)s\partial_s+s^2(-|\eta|^2)-n][(s\partial_s)^2-(n-1)s\partial_s+s^2(-|\eta|^2)+\frac{n^2-4}{2}],           
\end{align*}
and then
\begin{align*}
L_0(t,\,\widehat{\eta})\,&=\,[(t\partial_t)^2\,-\,(n-1)t\partial_t\,-\,t^2\,-\,n][(t\partial_t)^2\,-\,(n-1)t\partial_t\,-\,t^2\,+\,\frac{n^2-4}{2}]\\
&=\,L_1\,\circ\,L_2,
\end{align*}
with
$s\partial_s\,=\,s|\eta|\,\partial_{s|\eta|}\,=\,t\,\partial_t,$ and
$L_0$ here does not depend on $\tilde{y}$. Now we have used the full                                 
symmetry of the operator, and made it into the simplest form.

Let us consider the relationship of Fredholm property among $N(L)$,
$\widehat{N(L)}$ and $L_0$,  in $t^{\delta}L^2$, for
$\delta\,>\,\frac{n}{2}$. We know that the first two operators have
the same properties of injectivity and surjectivity. Let
\begin{align*}
L_0\varphi(t)=\,0,
\end{align*}
by definition, it holds if and only if
\begin{align*}
\widehat{N(L)}\varphi(s|\eta|)\,=\,0.
\end{align*}
But then
\begin{align*}
\widehat{N(L)}(a(\eta)\varphi(s\,|\eta|))\,=\,a(\eta)\,\widehat{N(L)}\varphi(s\,|\eta|)\,=\,0,
\end{align*}
for all $a(\eta)$ smooth, since the derivative is only in $s$
direction, with fixed $\eta$. Then,                                       
using the inverse Fourier transformation,
\begin{align*}
N(L)\,\int_{\mathbb{R}^{n-1}}e^{2\pi i\,\langle y,\,\eta
\rangle}\,a(\eta)\,\varphi(s|\eta|)\,d\,\eta\,=\,0.
\end{align*}
This means kernel of one dimensional $L_0$ corresponds to the             
infinite dimensional kernel of $N(L)$, and this construction also
gives the fact that the kernel of $N(L)$ is either trivial or of
infinite dimension. But if $\widehat{N(L)}$ is injective, then $L_0$
is injective. Conversely, if $L_0$ is injective, then
$\widehat{N(L)}$ is injective, and so is $N(L)$. We have a dual argument of the surjectivity for $\delta\,<\,\frac{n}{2}$. As in      
\cite{Mazzeo1}, $L_0$ is Fredholm when $\delta\,\notin \Lambda$,
with the set $\Lambda$ in (\ref{indicesset}), and $N(L)$ is                         
semi-Fredholm with either infinite dimensional kernel or cokernel.
Roughly speaking, $L$ is a small perturbation of $N(L)$ near $\partial M$. When $N(L)$                
is injective or surjective, $L$ is essentially injective or
essentially surjective, which will be Theorem \ref{SLSemiFredholm}                                             
and Theorem \ref{HLSemiFredholm}.

To see the semi-Fredholm property of $L$, the strategy is to first
study the Fredholm property of $L_0$ and $N(L)$, and finally obtain
the semi-Fredholm property of $L$ using Mazzeo's theorems which we
list here as Theorem \ref{SFPSL} and Corollary \ref{LHS}.

Now we discuss on the Fredholm property of $L_0$, $L_1$ and $L_2$ on
the weighted spaces. To this end, we introduce Bessel functions as
solutions to the Bessel equation as follows, which is well studied,
\begin{align*}
x^2\frac{d^2y}{dx^2}\,+\,x\frac{dy}{dx}\,-\,(x^2\,+\,\alpha^2)y\,=\,0,
\end{align*}
where $\alpha$ is a complex number.

The Bessel functions $I_{\alpha}$ and $I_{-\alpha}$ form a basis of
linear space of solutions to the Bessel function above, while
$\{I_{\alpha},\,K_{\alpha}\}$ is another basis.  For
$\text{Re}(\alpha)\,>\,-\,\frac{1}{2}$, and                         
$-\,\frac{\pi}{2}\,<\,arg(x)\,<\,\frac{\pi}{2}$, the integral
representations of these solutions are as follows,
\begin{align*}
I_{\alpha}(x)\,&=\,\frac{(\frac{x}{2})^{\alpha}}{\Gamma(\alpha+\frac{1}{2})\Gamma(\frac{1}{2})}\int_{-1}^1\,e^{-x\,t}(1\,-\,t^2)^{\alpha\,-\,\frac{1}{2}}\,dt,\\  
K_{\alpha}(x)\,&=\,\frac{\pi}{2}\,\frac{I_{\alpha}(x)\,-\,I_{-\alpha}(x)}{\sin(\alpha\,\pi)}\,=\,\frac{\Gamma(\frac{1}{2})(\frac{x}{2})^{\alpha}}{\Gamma(\alpha+\frac{1}{2})}\int_1^{\infty}\,e^{-x\,t}(t^2\,-\,1)^{\alpha\,-\,\frac{1}{2}}\,dt,
\end{align*}
with $x$ a complex number. See Page 172 and 77 in \cite{GW}. Note                                     
that $I_{\alpha}$ is bounded near $x=0,$ and it increases
exponentially near $+\infty$, and                                                                     
\begin{align*}
K_{\alpha}(x)\,\sim\,C(\varepsilon)x^{\text{Re}(\alpha)}e^{-\,x\,+\,\varepsilon},
\end{align*}
for any $\varepsilon\,>\,0$, as $x\,\rightarrow\,+\,\infty$. Also
$K_{\alpha}(x)$ is bounded for $\text{Re}(\alpha)\,\geq\,0$, near
$x\,=\,0$. The form $K_{\alpha}(x)$ is more useful near
$x\,=\,\infty$, since it decays exponentially.

We want to solve the following ODE, by transferring it into the
Bessel type equations as above.
\begin{align*}
L_1\,u\,=\,((t\partial_t)^2\,-\,(n-1)t\partial_t\,-\,t^2\,-\,n)\,u\,=\,0.
\end{align*}
Let $u\,=\,t^{\beta}\,\widetilde{u}$, then we obtain that
\begin{align}\label{Besselfunction}
t^{\beta}((t\partial_t)^2\,\widetilde{u}\,+\,(2\,\beta\,+\,1\,-\,n)\,t\partial_t\,\widetilde{u}\,+\,(-n\,-\,t^2\,\,+\,{\beta}^2\,-\,\beta(n-1))\widetilde{u})\,=\,0.
\end{align}
Then, letting $2\,\beta\,+\,1\,-\,n\,=\,0$, the equation
(\ref{Besselfunction}) is just the form of the Bessel function
defined as above. In this case, $\beta\,=\,\frac{n-1}{2}$, and then
the index $\alpha\,=\,\frac{n+1}{2}$.

\noindent Therefore,
\begin{align*}
u(t)\,=\,t^{\frac{n-1}{2}}(C_1I_{\frac{n+1}{2}}(t)\,+\,C_2K_{\frac{n+1}{2}}(t)).
\end{align*}
In fact,
\begin{align*}
t^{\frac{n-1}{2}}I_{\frac{n+1}{2}}(t|\eta|)\sim\,t^{n}|\eta|^{\frac{n+1}{2}},
\,\,\,\,\,
\,t^{\frac{n-1}{2}}K_{\frac{n+1}{2}}(t|\eta|)\sim\,t^{-1}|\eta|^{-\frac{n+1}{2}},
\end{align*}
near $t\,=\,0$. Moreover,
\begin{align*}
t^{\frac{n-1}{2}}I_{\frac{n+1}{2}}(t|\eta|)\sim\,t^{\frac{n}{2}-1}e^{t|\eta|}/\sqrt{2\pi
|\eta|},\,\,\,\,
t^{\frac{n-1}{2}}K_{\frac{n+1}{2}}(t|\eta|)\sim\,t^{\frac{n}{2}-1}e^{-t|\eta|}\sqrt{\frac{\pi}{2\,|\eta|}},
\end{align*}
as $t\,\to\,\infty.$ 

Similarly,
\begin{align*}
L_2\,u\,=\,((t\partial_t)^2\,-\,(n-1)t\partial_t\,-\,t^2\,+\,\frac{n^2-4}{2})\,u\,=\,0.
\end{align*}
Let $u(t)\,=\,t^{\beta}\,\widetilde{u}(t)$, then
\begin{align*}
t^{\beta}((t\partial_t)^2\,\widetilde{u}\,+\,(2\,\beta\,+\,1\,-\,n)\,t\partial_t\,\widetilde{u}\,+\,(\frac{n^2-4}{2}\,-\,t^2\,\,+\,{\beta}^2\,-\,\beta(n-1))\widetilde{u})\,=\,0.
\end{align*}
Set $2\,\beta\,+\,1\,-\,n\,=\,0$, so that $\beta\,=\,\frac{n-1}{2}$,
and then $\tilde{u}$ is a solution to the Bessel equation with
$\alpha\,=\,\frac{i\,\sqrt{n^2+2n-9}}{2}$.
\begin{align*}
u(t)\,=\,t^{\frac{n-1}{2}}\,(C_1I_{\frac{i\,\sqrt{n^2+2n-9}}{2}}(t)\,+\,C_2K_{\frac{i\,\sqrt{n^2+2n-9}}{2}}(t)).
\end{align*}

\noindent By the expansion of the series form of the Bessel
functions, as in [\cite{N.N.L}, P. 108], we have
\begin{align*}
t^{\frac{n-1}{2}}I_{\alpha}(t|\eta|)\sim\,t^{\frac{n-1}{2}+\alpha}|\eta|^{\alpha}/(2^{\alpha}\Gamma(1\,+\,\alpha)),
\end{align*}
and
\begin{align*}
t^{\frac{n\,-\,1}{2}}I_{-\,\alpha}(t|\eta|)\sim\,t^{\frac{n-1}{2}-\alpha}|\eta|^{\alpha}/(2^{\alpha}\Gamma(1\,-\,\alpha)),
\end{align*}
with $\alpha\,=\,\frac{i\,\sqrt{n^2+2n-9}}{2}$, near $t\,=\,0$. Now                    
it is easy to see that the linear
 combination
 \begin{align*}
 x^{\frac{n\,-\,1}{2}}(C_1\,x^{i\,\frac{\sqrt{n^2\,+\,2n\,-\,9}}{2}}\,+\,C_2\,x^{-\,i\,\frac{\sqrt{n^2\,+\,2n\,-\,9}}{2}})                
 \end{align*}
 can never vanish to infinite order at $t\,=\,0$ if either $C_1 \neq 0$ or $C_2 \neq 0$. Also,                                             
\begin{align*}
t^{\frac{n-1}{2}}K_{\alpha}(t|\eta|)\sim\,t^{\frac{n-1}{2}}\,\frac{\pi}{2}\frac{I_{\alpha}(t|\eta|)\,-\,I_{-\,\alpha}(t|\eta|)}{\sin(\alpha\pi)},
\end{align*}
with $\alpha\,=\,\frac{i\,\sqrt{n^2+2n-9}}{2}$, and $|\eta| \neq 0$,                                                                       
near $t\,=\,0$.

Using the integral form as above, we have that
$I_{\alpha}(t)$ grows exponentially, while $K_{\alpha}$ decays
exponentially as 
$t\,\rightarrow\,+\,\infty$, for
$\alpha\,=\,\frac{i\,\sqrt{n^2+2n-9}}{2}$.  

Denote $L_0^t$ to be the $L^2$ adjoint of $L_0$ in the measure
$d\,t$, and
\begin{align*}
L_0^{*}\,=\,t^{2\delta}L_0t^{-2\delta},
\end{align*}
to be the adjoint of $L_0$ in $t^{\delta}L^2$ in the measure
$t^{-2\delta}\,d\,t$. These are all elliptic operators, with
boundary spectra:
\begin{align*}
&spec_b(L_0^t)\,=\,\{\,-\zeta\,-\,1\,:\,\zeta\,\in\,spec_b(L_0)\},\\
&spec_b(L_0^{*})\,=\,\{\,-\zeta\,+\,2\delta\,-\,1\,:\,\zeta\,\in\,spec_b(L_0)\}.
\end{align*}
For example, for                                                                                    
$L_1\,=\,(t\partial_t)^2\,-\,(n-1)(t\partial_t)\,-\,t^2\,-\,n$,
\begin{align*}
\int \,L_1u\,v\,dt\,=\,\int\,u\,L_1^tv\,dt.
\end{align*}
Then
\begin{align*}
L_1^t\,=\,(-\partial_t\,(t\,\cdot))^2\,+\,(n-1)(\partial_t\,(t\,\cdot))\,-\,t^2\,-\,n,
\end{align*}
with
\begin{align*}
\partial_t\,(t\,\cdot)\,=\,t\partial_t\,+\,1,
\end{align*}
and $p^t(\xi)\,=\,p(-(\xi\,+\,1))$, for the quadratic polynomial
$p$. Also, for $L_1^{*}$, using the fact that
\begin{align*}
-\partial_t(\,t\,t^{-2\delta}\cdot)\,=\,-\,t^{-2\delta}(-2\delta\,+\,1\,+\,t\,\partial_t)\,=\,t^{-2\delta}(2\delta\,-\,1\,-\,t\,\partial_t),
\end{align*}
and
\begin{align*}
\int\,L_1\,u\,v\,t^{2\delta}\,dt\,=\,-\,\int\,u\,t^{-2\delta}L_1^{t}(t^{2\delta}v)\,t^{2\delta}\,dt,
\end{align*}
we obtain the boundary spectra as listed above. For the fourth order
differential equation, we have obtained four linearly independent
solutions, and they generalize the solution space.

Let $\delta\,=\,\frac{n-1}{2}\,+\,\frac{1}{2}=\frac{n}{2}$, we have                       
$L_1^{*}\,=\,L_1$, and $L_2^{*}\,=\,L_2$.








\begin{defn}
We say that an operator $L$ has the unique continuation property on
a boundary $B$ if any solution of $L\,u\,=\,0$ vanishing to infinite
order at $B$ vanishes identically.
\end{defn}

\begin{hyp}\label{hypr}                                  
For each $\tilde{y}$ and $\hat{\eta}$, both $L_0$ and its adjoint
$L_0^{*}$ (The dual of $L_0$ with respect to the space
$t^{\text{Re}\delta}L^2$ for any $\delta$ we need) have the unique
continuation property at $\{t\,=\,0\}$.
\end{hyp}

We know from the discussion above that $L_0$ satisfies the unique
continuation property. Under the continuation hypothesis, we have
that for each element $(\tilde{y},\,\hat{\eta})\,\in\,N_0$, $L_0$ is
surjective on $x^{\delta}L^2$ or injective on $x^{\delta}L^2$ when                                                
$\delta$ is sufficiently negative or sufficiently large. For our
case, we use $\delta\,=\,\frac{n}{2}$ in Hypothesis \ref{hypr}. Now                                               
let us define $\overline{\delta}$ to be the minimal value of
$\delta$ so that $L_0$ is injective, and meanwhile
$\underline{\delta}$ the maximal value so that $L_0$ is surjective
dually. These values must lie in $\Lambda$. The following theorem
and corollary tell us the relationship between semi-Fredholm
properties of $L$ and the Fredholm properties of $L_0$, for certain
cases we need.

\begin{thm}\label{SFPSL}
$($\text{Theorem} $6.1.$ \text{in} \cite{Mazzeo1}$)$ Suppose
$L\,\in\,\text{Diff}_e^{m}(M)$ is elliptic and satisfies the unique
continuation hypothesis, and that $spec_b(L)$ is discrete. Suppose
also that $\delta\,\notin\,\Lambda$ is chosen so that either                                                  
$\delta\,>\,\bar{\delta}$ or $\delta\,<\,\underline{\delta}$. Then
$\,L:\,x^{\delta}H_e^{r+m}(M)\,\rightarrow\,x^{\delta}H_e^{r}(M)$
has closed range, and it is either essentially surjective, or
essentially injective, which means respectively that $L$ has either
an at most finite dimensional nullspace, or a finite dimensional                                             
cokernel. Therefore, it admits a generalized inverse $G$ and
orthogonal projectors $P_i$ onto the nullspace and orthogonal
complement of the range of $L$ which are edge operators, such that,
\begin{align*}
GL\,=\,I\,-\,P_1,\\
LG\,=\,I\,-\,P_2.
\end{align*}

\end{thm}
Since the edge operators used in the proof of the weighted Sobolev
spaces are bounded in the appropriate H$\ddot{\text{o}}$lder spaces,
the corresponding result for H$\ddot{\text{o}}$lder spaces follows.
\begin{cor}\label{LHS}
$($\text{Corollary} $6.4.$ \text{in} \cite{Mazzeo1}$)$ For $L$ as in
Theorem $\ref{SFPSL}$, $k\,\geq\,m$ a positive integer and
$0\,<\,\alpha\,<1$ the mapping
$L\,:\,x^{\nu}\Lambda^{k,\,\alpha}\,\rightarrow\,x^{\nu}\Lambda^{k\,-\,m,\,\alpha}$                          
is semi-Fredholm provided $\nu\,=\,\delta\,-\,\frac{1}{2}$ and
$\delta\,\notin\,\Lambda$ is as in the previous theorem. If
$\delta\,<\,\underline{\delta}$ or $\delta\,>\,\bar{\delta}$ so that
$L$ is essentially surjective or essentially injective, then
topologically, we have the splitting,
\begin{align*}
&x^{\nu}\Lambda^{k,\,\alpha}\,=\,P_1(x^{\nu}\Lambda^{k,\,\alpha})\,\oplus\,(I\,-\,P_1)(x^{\nu}\Lambda^{k,\,\alpha}),\\           
&x^{\nu}\Lambda^{k\,-\,m,\,\alpha}\,=\,P_2(x^{\nu}\Lambda^{k\,-\,m,\,\alpha})\,\oplus\,(I\,-\,P_2)(x^{\nu}\Lambda^{k\,-\,m,\,\alpha}).
\end{align*}
\end{cor}
Let us compute $\overline{\delta}$ and $\underline{\delta}$ for
$L_0$. First, for $L_1$, since
$t^{\frac{n-1}{2}}I_{\frac{n+1}{2}}(t|\eta|)$ increases
exponentially as $t$ goes to $\infty$ (here $|\eta|\,\neq\,0$), it
does not lie in $t^{\delta}L^2$ for any $\delta\,>\,0$; furthermore,
\begin{align*}
t^{\frac{n-1}{2}}K_{\frac{n+1}{2}}(t|\eta|)\,\in\,t^{\delta}L^2(\mathbb{R}_{+}),
\end{align*}
for $\delta\,<\,-\,\frac{1}{2}$. Similarly, for $L_2$,
$t^{\frac{n-1}{2}}I_{\frac{i\,\sqrt{n^2+2n-9}}{2}}(t|\eta|)$ grows
exponentially when $t$ goes to $\infty$ ( with $|\eta|\,\neq\,0$),
and
\begin{align*}
t^{\frac{n-1}{2}}K_{\frac{i\,\sqrt{n^2+2n-9}}{2}}(t|\eta|)\,\in\,t^{\delta}L^2(\mathbb{R}_{+}),
\end{align*}
for $\delta\,<\,\frac{n-1}{2}+\frac{1}{2}\,=\,\frac{n}{2}$.
Therefore, $L_1$ and $L_2$ both have trivial kernel in the space
$x^{\delta}L^2(M,\,\sqrt{dxdy})$ for $\delta\,>\,\frac{n}{2}$. But
$\text{Ker}(L_2)$ is nontrivial for $\delta\,<\,\frac{n}{2}$. Also
the composition of two injective map is still injective. Therefore,
$\overline{\delta}\,=\,\frac{n}{2}$ for $L_0\,=\,L_1\,\circ\,L_2$.
Since $L_0$ is self-adjoint in $t^{\frac{n}{2}}L^2(\mathbb{R}_{+})$,
we have that $\underline{\delta}\,=\,\frac{n}{2}$. Since it
satisfies the conditions of Theorem \ref{SFPSL} and Corollary
\ref{LHS}, therefore Theorems \ref{SLSemiFredholm} and
\ref{HLSemiFredholm} are proved. \qed

To conclude this section, we want to see when $L$ is injective or                                             
surjective in the
special case of Poincar$\acute{\text{e}}$-Einstein manifolds. For a                                          
Poincar$\acute{\text{e}}$-Einstein manifold $(M,\,g)$ with
$g\,=\,x^{-2}h$, without loss of                                                                                 
generality we assume $R_g\,=\,-n(n-1)$. Let us first consider it in
the weighted Sobolev spaces. We have
\begin{align*}
L\,=\,(\Delta_g\,-\,n)(\Delta_g\,+\,\frac{(n+2)(n-2)}{2})\,=\,\mathcal
{T}_1\,\circ\,\mathcal {T}_2.
\end{align*}
We know that $L$ is self-adjoint with respect to
$x^{\frac{n}{2}}L^2(M,\,\sqrt{dx\,dy})$. Then to show that
$L\,:\,x^{\delta}H_e^4(M)\,\rightarrow\,x^{\delta}L^2(M)$ is
surjective for $0\,\leq\,\delta\,<\,\frac{n}{2}$, we only need to                             
show that $L$ is injective when $\delta\,>\,\frac{n}{2}$. For that,
we only need to show that $\mathcal {T}_1$ and $\mathcal {T}_2$ are
injective for $\delta\,>\,\frac{n}{2}$.

As a special case, if $\mathcal {T}_1$ and $\mathcal {T}_2$ are
injective in $L^2(M,\,g)$, which is
$x^{\frac{n}{2}}L^2(M,\,\sqrt{dxdy})$, then we are done. For the                               
Poincar$\acute{\text{e}}$ ball $(B,\,g_{-1})$, we know that the
Laplacian $-\,\Delta_g$ has pure continuous spectrum, consisting of
$[\frac{(n-1)^2}{4},\,\infty)$, with
$\lambda_0\,=\,\frac{(n-1)^2}{4}$. So $\mathcal {T}_1$ is injective
when $\delta\,>\,\frac{n}{2}$ in this special case.

As mentioned in the introduction, the way we prove the surjectivity
of                                                       
 $\mathcal {T}_2$ in the following involves unique continuation property and boundary smoothness argument,
 which need $x$ and $h$ to be smooth enough.               
So we assume 
the defining function $x$ and the metric $h$ to be smooth up to                          
the boundary.

\begin{lem}\label{injectivitylem}
$\mathcal {T}_1,\,\mathcal {T}_2
:\,x^{\delta}H_e^{2+m}(M,\,\sqrt{dxdy})\,\rightarrow\,x^{\delta}H_e^{m}(M,\,\sqrt{dxdy})$,
are both injective for $\delta\,>\,\frac{n}{2}$, and all
$m\,\geq\,0$.
\end{lem}
\noindent Proof of Lemma \ref{injectivitylem}. By the regularity
argument, we only need to discuss on the case $m\,=\,0$. The proof
is as follows.

If $u\,\in\,x^{\delta}L^2(M,\,\sqrt{dxdy})$, for
$\delta\,>\,\frac{n}{2}$, then $u\,\in\,L^2(M,\,g)$. Moreover, if                                 
also
\begin{align*}
\mathcal {T}_1\,u\,=\,(\Delta_g\,-\,n)\,u\,=\,0,
\end{align*}
by Weyl's lemma, $u\,\in\,H^1(M,\,g)$. Now we multiply $u$ on both
sides of the equation, and integrate by parts, and then we have
\begin{align*}
-\int_M\,(|\nabla u|_g^2\,+\,u^2)\,dV_g\,=\,0.
\end{align*}
Therefore, $u\,=\,0$. Then we have that $\mathcal {T}_1$ is
injective.

Let us consider the equation
\begin{align}\label{linearequation2}
\mathcal {T}_2\,u\,=\,(\Delta_g\,+\,\frac{n^2\,-\,4}{2})\,u\,=\,0.                                                            
\end{align}
First, we know from \cite{Mazzeo} that $(-\Delta_g\,-\,\lambda)$                   
satisfies the unique continuation property, for any constant real
number $\lambda$(See Corollary 11 in \cite{Mazzeo}), namely, if                    
$u\,\in\,C^{\infty}(M)$ satisfies
$(-\,\Delta_g\,-\,\lambda)\,u\,=\,0$, and $u$ vanishes to infinite
order along an open set of $\partial M$, then $u\,=\,0$. Moreover,                 
in \cite{Mazzeo}, combining the boundary regularity result and the
unique continuation result for $(-\Delta_g\,-\,\lambda)$, it                       
was proved in \cite{Mazzeo} that if $\lambda\,>\,\frac{(n+1)^2}{4}$, 
$u\,\in\,L^2(M,\,g)$ and $(-\Delta_g\,-\,\lambda)\,u\,=\,0$ then
$u\,=\,0$. It is easy to check that when $n\,\geq\,5$,
$\frac{n^2\,-\,4}{2}\,>\,\frac{(n+1)^2}{4}$. Therefore, for
$n\,\geq\,5$, $\mathcal {T}_2$ is injective in
$L^2(M,\,g)\,=\,x^{\frac{n}{2}}L^2(M,\,\sqrt{dxdy})$.

When $n\,=\,4$, since $\frac{n^2\,-\,4}{2}\,<\,\frac{(n+1)^2}{4}$,                           
we can not use his result directly. But since we still have the unique                       
continuation property for $\mathcal {T}_2$, we only need to prove
the boundary regularity of $u$. Actually, for our case we do not                             
require $u\,\in\,L^2(M,\,g)$ but allow
$u\,\in\,x^{\delta}L^2(M,\,\sqrt{dx\,dy})$, for
$\delta\,>\,\frac{n}{2}$, and the method of proving the boundary
regularity still works here. For completeness, we give the details
here. Assume $u\,\in\,x^{\delta}L^2(M,\,\sqrt{dxdy})$ for
$\delta\,>\,\frac{n}{2}$ satisfies equation (\ref{linearequation2}).
Let us denote the indicial roots of $\mathcal {T}_2$ as $s_1,\,s_2$.
Using the boundary expansion from Section 7 in \cite{Mazzeo1}, we
have that
\begin{align}\label{generalformalexpansion}                                                 
u\,\sim\,\displaystyle\sum_{j=0}^{\infty}\,(\sum_{p\,=\,0}^{N_1}x^{s_1+j}(\log(x))^pu_{1,j,p}(y)\,+\,\sum_{p\,=\,0}^{N_2}x^{s_2+j}(\log(x))^pu_{2,j,p}(y)).
\end{align}
Since the real parts of the indicial roots $s_1,\,s_2$ are both                               
$\frac{n-1}{2}$, which is less than $\delta\,-\,\frac{1}{2}$, then
by Theorem (7.17) in \cite{Mazzeo1}, we have that $u_{1,\,0,\,p}$
and $u_{2,\,0,\,p}$ vanish for all $p$, and that the coefficients                              
$u_{i,j,p}(y)$ are all smooth, and by Theorem (7.3) in
\cite{Mazzeo1}, using the substitution of the expansion into the
equation, we have that the coefficients all vanish by induction.
Then, by the unique continuation property, we have that $u\,=\,0$.

It follows that $\mathcal {T}_1$ and $\mathcal {T}_2$ are both
injective when $\delta\,>\,\frac{n}{2}$. This proves the lemma. \qed
\vskip2pt

The lemma implies that $L$ is injective for $\delta\,>\,\frac{n}{2}$
on the Poincar$\acute{\text{e}}$-Einstein manifolds. Since $L$ is
self-adjoint in $x^{\frac{n}{2}}L^2(M,\,\sqrt{dx\,dy})$, then $L$ is
surjective when $0<\,\delta\,<\,\frac{n}{2}$.

The linear edge operators used above are all bounded linear
operators in the weighted H$\ddot{\text{o}}$lder spaces, and can be
used correspondingly in the weighted H$\ddot{\text{o}}$lder spaces.
Then the corresponding statement for the weighted                              
H$\ddot{\text{o}}$lder spaces is as follows. Let
\begin{align*}
L\,:\,x^{\nu}\Lambda^{4,\,\alpha}(M)\,\rightarrow\,x^{\nu}\Lambda^{0,
\,\alpha}(M).
\end{align*}
Here $0<\alpha<1$. Then $L$ is injective when                                                                                       
$\nu\,=\,\delta\,-\,\frac{1}{2}\,>\,\frac{n-1}{2}$, while $L$ is
surjective when
$0\,<\,\nu\,=\,\delta\,-\,\frac{1}{2}\,<\,\frac{n-1}{2}$ on the                
Poincar$\acute{\text{e}}$-Einstein manifolds $M$.
\begin{Remark}\label{Remark_Regularity}
Generally, on an asymptotically hyperbolic manifold $(M,\,g)$ with a smooth defining function $x$ and $h\,=\,x^2g$ smooth, let $u\,\in\,\text{Ker}(L)$ for $L$ defined in (\ref{linearoperator}) in the weighted         
H$\ddot{\text{o}}$lder spaces
$x^{\nu}\Lambda^{4,\,\alpha}(M,\,\sqrt{dx\,dy})$, for
$0\,<\,\nu\,<\,\frac{n-1}{2}$ and $0<\alpha<1$. Then,                                                                               
$u\,\in\,x^{\nu}\Lambda^{m,\,\alpha}$ for all $m\in\,\mathbb{N}$,
and $u$ has the following weak expansion with                                                                                       
coefficients which are generally distributions,                                                                                           
\begin{align}\label{sequence0}
u(x,\,y)\,\sim\,\displaystyle\sum_{j\,=\,0}^{+\,\infty}(u_{0j}(y)x^{\frac{n-1}{2}\,+\,i\frac{\sqrt{n^2+2n-9}}{2}\,+\,j}\,+\,u_{1j}(y)x^{\frac{n-1}{2}\,-\,i\frac{\sqrt{n^2+2n-9}}{2}\,+\,j}\,+\,x^{n\,+\,j}u_{2j}(y)),                 
\end{align}
in the sense that
\begin{align*}
u(x,y)\,-\,\displaystyle\sum_{j=0}^k(u_{0j}(y)x^{\frac{n-1}{2}\,+\,i\frac{\sqrt{n^2+2n-9}}{2}\,+\,j}\,+\,u_{1j}(y)x^{\frac{n-1}{2}\,-\,i\frac{\sqrt{n^2+2n-9}}{2}\,+\,j})\,=\,o(x^{\frac{n-1}{2}+k}),
\end{align*}
for $k\,\geq\,0$. If either $u_{00}$ or                                                  
$u_{01}$ is smooth, then all the coefficients are smooth. The more
precise regularity of the coefficients in a weighted Sobolev space
setting can be found in Chapter $7$ in \cite{Mazzeo1}.                                   
\end{Remark}
\begin{Remark}\label{Remark_smoothness}
On a Poincar$\acute{\text{e}}$-Einstein manifold $(M,\,g)$ with a smooth defining function $x$ and $h\,=\,x^2g$ smooth, for $0\,<\,\nu\,<\frac{n-1}{2}$ and $0<\alpha<1$, since                
$\mathcal {T}_1$ is injective, an element $u$ in the kernel of $L$
is exactly an element in the kernel of $\mathcal {T}_2$. By
Proposition $3.4.$ in \cite{Graham and Zworski}, for any chosen
$u_{00}\,\in\,C^{\infty}$ or $u_{10}\,\in\,C^{\infty}$, there exists                             
a unique $u\,\in\,x^{\nu}\Lambda^{4,\,\alpha}(M,\,\sqrt{dx\,dy})$,
for $0\,<\,\nu\,<\,\frac{n-1}{2}$, in the kernel of $L$, so that $u$
has the expansion $(\ref{sequence0})$ with smooth coefficients.
\end{Remark}
\vskip5pt

\section {The nonlinear problem}
\vskip10pt

Now let us return to the perturbation problem. It is more convenient
to work in weighted H$\ddot{\text{o}}$lder spaces. 
Let $(M,\,g)$ be an asymptotically hyperbolic manifold defined as in                         
the introduction. Let $\tilde{g},\,u$ also be defined as in the
introduction, and let the prescribed curvature
$Q_{\tilde{g}}\,=\,f$. 
Define the operator $\mathcal
{T}\,:\,x^{\nu}\Lambda^{4,\,\alpha}(M)\,\rightarrow\,x^{\nu}\Lambda^{0,\,\alpha}(M)$
as follows,

\begin{align*}
\mathcal {T}(u)\,=\,\left\{\begin{array}{l@{\quad\quad}l} 2\,f\,e^{4u}\,-\,2\,Q_g\,-\,8\,Q_g\,u,\,\, & n\,=\, 4,\\ \\
 \frac{n-4}{2}\,(1\,+\,u)^{\frac{n+4}{n-4}}\,f\,-\frac{n-4}{2}Q_g\,-\,\frac{n+4}{2}Q_g\,u, & n\,\geq\,5.\end{array} \right.
\end{align*}
We rewrite it in the form
\begin{align*}
\mathcal {T}(u)=\left\{\begin{array}{l@{\quad\quad}l}
2(e^{4u}-1-4u)f+2(f-Q_g)-8(Q_g-f)\,u, & n\,=\, 4,\\ \\                                                                      
 \frac{n-4}{2}((1+u)^{\frac{n+4}{n-4}}-1-\frac{n+4}{n-4}u)f+\frac{n-4}{2}(f-Q_g)+\frac{n+4}{2}(f-Q_g)u,& n\,\geq\,5.\end{array} \right.
\end{align*}
Let $L$ be as in (\ref{linearoperator}), then the prescribed                                                                
$Q$-Curvature equation is
\begin{align}\label{general_version_nonlinear_equation}
L\,u\,=\,\mathcal {T}(u).
\end{align}
Let $0\,<\,\nu\,<\,\underline{\nu}\,=\,\frac{n-1}{2}$ and $0<\alpha<1$, so that $L$                             
is essentially surjective. Moreover, in the following we assume that
$L$ is surjective. Then
\begin{align*}
L:\,V_1\,=\,(I\,-\,P_1)(x^{\nu}\Lambda^{4,\,\alpha}(M))\,\to\,x^{\nu}\Lambda^{0,\,\alpha}(M)
\end{align*}
is an isomorphism, using topological splitting of
$x^{\nu}\Lambda^{4,\,\alpha}(M)$ in Theorem \ref{HLSemiFredholm} and                                                           
the open mapping theorem. That is,
\begin{align}\label{uniformlybounded}
C_1\,\|u\|_{x^{\nu}\Lambda^{4,\alpha}(M)}\,\leq\,\|L\,u\|_{x^{\nu}\Lambda^{0,\alpha}(M)}\,\leq\,C_2\,\|u\|_{x^{\nu}\Lambda^{4,\alpha}(M)},
\end{align}
for some constant $C_2\,>\,C_1\,>\,0$, for all $u\,\in\,V_1$. We                                                              
denote the inverse of $L$ as
\begin{align*}
L^{-1}:\,x^{\nu}\Lambda^{0,\,\alpha}(M)\,\rightarrow\,V_1.
\end{align*}
Let $f\,\in\,C^{\alpha}(M)$, 
and
\begin{align*}
(Q_g\,-\,f)\,\in\,x^{\nu}\Lambda^{0,\,\alpha},                                                           
\end{align*}
with its small norm to be determined later.                                                    
We want to use elements in kernel of $L$ to parametrize the
perturbation solutions to the nonlinear problem at $0$. We will
define a new map for each element in the kernel of $L$, and use it
to construct a contraction map. For any fixed
$u_1\,\in\,\text{Ker}(L)$, for any $u_2\,\in\,V_1$, let
$u\,=\,u_1\,+\,u_2$, and
\begin{align*}
\mathcal{T}_{u_1}(u_2)\,=\,\mathcal {T}(u_1\,+\,u_2).
\end{align*}
Now $L^{-1}\,\circ\,\mathcal {T}_{u_1}\,:\,V_1\,\rightarrow\,V_1$.

From now on, let $u_1$ be any fixed element in
$B_{\epsilon}(0)\,\bigcap\,\text{Ker}(L)$, and
$u_2\,\in\,B_{\epsilon}(0)\,\bigcap\,V_1$, with small
$\epsilon\,\in\,(0,\,1)$ to be determined. Note that
\begin{align*}
\|\mathcal
{T}_{u_1}(u_2)\|_{x^{\nu}\Lambda^{0,\alpha}}\leq\left\{\begin{array}{l@{\quad\quad}l}2\|(e^{4u}-1-4u)f\|_{x^{\nu}\Lambda^{0,\alpha}(M)}+2\|(f-Q_g)\|_{x^{\nu}\Lambda^{0,\alpha}(M)}\\                
+8\|(f-Q_g)u\|_{x^{\nu}\Lambda^{0,\alpha}(M)}, & n=4,\\ \\
 \frac{n-4}{2}\|((1+u)^{\frac{n+4}{n-4}}-1-\frac{n+4}{n-4}u)f\|_{x^{\nu}\Lambda^{0,\alpha}(M)}\\
 +\frac{n-4}{2}\|(f-Q_g)\|_{x^{\nu}\Lambda^{0,\alpha}(M)}+\frac{n+4}{2}\|(f-Q_g)u\|_{x^{\nu}\Lambda^{0,\alpha}(M)},& n\geq
5.\end{array} \right.
\end{align*}
Then we have
\begin{align*}
\|\mathcal
{T}_{u_1}(u_2)\|_{x^{\nu}\Lambda^{0,\,\alpha}}\,\leq&\,C(n)\big(\|f\|_{L^{\infty}}\|(u_1\,+\,u_2)^2\|_{x^{\nu}\Lambda^{0,\,\alpha}}\,+\,(1+\|u_1+u_2\|_{L^{\infty}})\|f\,-\,Q_g\|_{x^{\nu}\Lambda^{0,\,\alpha}}\\
&+\,\|x^{-\nu}(u_1\,+\,u_2)\|_{L^{\infty}}\|(u_1\,+\,u_2)\|_{L^{\infty}}(\|f\|_{\Lambda^{0,\,\alpha}}\,+\,\|Q_g\|_{\Lambda^{0,\,\alpha}})\\
&+\,\|f\,-\,Q_g\|_{L^{\infty}}\|u_1\,+\,u_2\|_{x^{\nu}\Lambda^{0,\,\alpha}}\big).
\end{align*}
where $C\,>\,0$ is a constant depending only on $n$, the diameter of
$M$ and $\nu$. By the definition of the weighted norm,
\begin{align}\label{normcontrolling}
\|\phi\|_{L^{\infty}}\,\leq\,\|\phi\|_{\Lambda^{0,\,\alpha}},\,\,\text{and}\,\,\|\phi\|_{L^{\infty}}\,\leq\,C_0\,\|\phi\|_{x^{\nu}\Lambda^{0,\,\alpha}},
\end{align}
for a constant $C_0\,>\,0$ depending on the defining function and
$\nu$, for any $\phi\,\in\,x^{\nu}\Lambda^{0,\,\alpha}$. Therefore,
\begin{align*}
\|\mathcal{T}_{u_1}(u_2)\|_{x^{\nu}\Lambda^{0,\,\alpha}}\,\leq&\,C_1\,\big((\epsilon(\|f\|_{\Lambda^{0,\,\alpha}}\,+\,\|Q_g\|_{\Lambda^{0,\,\alpha}})+\|f-   
Q_g\|_{L^{\infty}})\|u_1\,+\,u_2\|_{x^{\nu}\Lambda^{0,\,\alpha}}\\
&+\,(1\,+\,\epsilon)\,\|f\,-\,Q_g\|_{x^{\nu}\Lambda^{0,\,\alpha}}\big),
\end{align*}
where $C_1$ depends on $n$, the defining function, the diameter of
$M$ and $\nu$, so that
\begin{align*}
\|L^{-1}\,\circ\,\mathcal
{T}_{u_1}(u_2)\|_{x^{\nu}\Lambda^{4,\,\alpha}}\,\leq&\,C\,\big((\epsilon(\|f\|_{\Lambda^{0,\,\alpha}}\,+\,\|Q_g\|_{\Lambda^{0,\,\alpha}})+\|f-Q_g\|_{L^{\infty}})\|u_1\,+\,u_2\|_{x^{\nu}\Lambda^{0,\,\alpha}}\\              
&+\,(1\,+\,\epsilon)\,\|f\,-\,Q_g\|_{x^{\nu}\Lambda^{0,\,\alpha}}\big),
\end{align*}
where $C\,=\,C_1\,\|L^{-1}\|$ depends on the defining function, the
diameter of $M$, $\nu$, $n$ and $\|L^{-1}\|$. We choose
$\epsilon\,\in\,(0,\,1)$ small so that
\begin{align}\label{condition_ineqn1}
16\,C\,\epsilon\,\|Q_g\|_{\Lambda^{0,\,\alpha}}\,<\,1,                         
\end{align}
and let $f$ satisfy that
\begin{align}\label{closequantity}
\|f\|_{\Lambda^{0,\,\alpha}}\,\leq\,2\,\|Q_g\|_{\Lambda^{0,\,\alpha}},\,\,\text{and}\,\,                                            
\|f\,-\,Q_g\|_{x^{\nu}\Lambda^{0,\,\alpha}}\,\leq\,\text{min}\{\frac{1}{4(1\,+\,\epsilon)C}\epsilon,\,\frac{\epsilon\|Q_g\|_{\Lambda^{0,\,\alpha}}}{C_0}\}.
\end{align}
Combining $(\ref{normcontrolling})$, we have                                    
\begin{align*}
\|L^{-1}\,\circ\,\mathcal
{T}_{u_1}(u_2)\|_{x^{\nu}\Lambda^{4,\,\alpha}}\,\leq\,\frac{3}{4}\epsilon.
\end{align*}
Therefore, $L^{-1}\,\circ\,{T}_{u_1}$ maps $B_{\epsilon}(0)\bigcap
V_1$ into $B_{\epsilon}(0)\bigcap V_1$.

For $u_3,\,u_4\,\in\,V_1\,\bigcap\,B_{\epsilon}(0)$,
\begin{align*}
&\|L^{-1}\,\circ\,\mathcal
{T}_{u_1}(u_3)\,-\,L^{-1}\,\circ\,\mathcal
{T}_{u_1}(u_4)\|_{x^{\nu}\Lambda^{4,\,\alpha}}\\
\leq&\,\|L^{-1}\|\,\|\mathcal {T}_{u_1}(u_3)\,-\,\mathcal              
{T}_{u_1}(u_4)\|_{x^{\nu}\Lambda^{0,\,\alpha}}\\
=&\left\{\begin{array}{l@{\quad\quad}l}\|L^{-1}\|\|2f(e^{4u_1}(e^{4u_3}-e^{4u_4})-4(u_3-u_4))-8(Q_g-f)(u_3-u_4)\|_{x^{\nu}\Lambda^{0,\alpha}},
& n=4,\\ \\   
 \|L^{-1}\|\|\frac{n-4}{2}((1+u_1+u_3)^{\frac{n+4}{n-4}}-(1+u_1+u_4)^{\frac{n+4}{n-4}}-\frac{n+4}{n-4}(u_3-u_4))f\\
 +\frac{n+4}{2}(f-Q_g)(u_3-u_4)\|_{x^{\nu}\Lambda^{0,\alpha}},&
n\,\geq\,5.\end{array}
 \right.
\end{align*}
But
\begin{align*}
e^{4(u_1\,+\,u_3)}\,-\,e^{4(u_1\,+\,u_4)}\,-\,4(u_3\,-\,u_4)\,=\,4\,(u_3\,-\,u_4)w,
\end{align*}
with
\begin{align*}
w=(\frac{e^{4(u_1+u_3)}-e^{4(u_1+u_4)}}{4(u_3-u_4)}-1)=(\int_0^1e^{4(u_1+u_4+t(u_3-u_4))}dt-1)\in
x^{\nu}\Lambda^{0,\alpha}\bigcap B_{C\epsilon}(0),                                                           
\end{align*}
with $C$ which does not depend on $u_3$, $u_4$, or
$\epsilon\,\in\,(0,\,1)$. We have similar results for $n\,\geq\,5$.                         
By the discussion above,
\begin{align*}
&\|L^{-1}\,\circ\,\mathcal
{T}_{u_1}(u_3)\,-\,L^{-1}\,\circ\,\mathcal
{T}_{u_1}(u_4)\|_{x^{\nu}\Lambda^{4,\,\alpha}}\\
&\leq\,\|L^{-1}\|\,\widetilde{C_0}\,(\,\epsilon\,\|f\|_{\Lambda^{0,\,\alpha}}\,\|u_3\,-\,u_4\|_{x^{\nu}\Lambda^{0,\,\alpha}}\,+\,\|Q_g\,-\,f\|_{x^{\nu}\Lambda^{0,\,\alpha}}\,\|u_3\,-\,u_4\|_{x^{\nu}\Lambda^{0,\,\alpha}})\\
&=\,\|L^{-1}\|\,\widetilde{C_0}\,(\,\epsilon\,\|f\|_{\Lambda^{0,\,\alpha}}\,+\,\|Q_g\,-\,f\|_{x^{\nu}\Lambda^{0,\,\alpha}})\,\|u_3\,-\,u_4\|_{x^{\nu}\Lambda^{0,\,\alpha}},
\,\,n\,\geq\, 4,
\end{align*}
where $\widetilde{C_0}$ depends only on the defining function, the                                                  
diameter of $M$, $\nu$ and $n$. Let $\epsilon$ be small so that
\begin{align}\label{condition_ineqn2}
8\widetilde{C_0}\|L^{-1}\|\,(1\,+\,\|Q_g\|_{\Lambda^{0,\,\alpha}})\,\epsilon\,<1,                                   
\end{align}
and let
\begin{align}\label{smallness}
\|Q_g\,-\,f\|_{x^{\nu}\Lambda^{0,\,\alpha}}\,\leq\,\frac{1}{8\widetilde{C_0}\|L^{-1}\|},
\end{align}
then we have
\begin{align*}
|L^{-1}\,\circ\,\mathcal {T}_{u_1}(u_3)\,-\,L^{-1}\,\circ\,\mathcal
{T}_{u_1}(u_4)\|_{x^{\nu}\Lambda^{4,\,\alpha}}\,\leq\,\frac{3}{8}\|u_3\,-\,u_4\|_{x^{\nu}\Lambda^{0,\,\alpha}}\\       
\leq\,\frac{3}{8}\|u_3\,-\,u_4\|_{x^{\nu}\Lambda^{4,\,\alpha}}.                                                        
\end{align*}
Note that $\|L^{-1}\|$ depends on the projection map $P_1$ that we
construct in Theorem \ref{HLSemiFredholm}. Therefore, if $L$ is
surjective for $\nu\,<\,\frac{n-1}{2}$, and also $\epsilon$ and $f$                      
satisfy the above conditions, then for each
$u_1\,\in\,B_{\epsilon}(0)\,\bigcap\,\text{Ker}(L)$,
\begin{align*}
L^{-1}\,\circ\,\mathcal
{T}_{u_1}\,:\,V_1\,\bigcap\,B_{\epsilon}(0)\,\rightarrow\,V_1\,\bigcap\,B_{\epsilon}(0)
\end{align*}
is a contraction map. This implies that there exists a unique
$u_2\,\in\,B_{\epsilon}(0)\,\bigcap\,V_1$, solving the equation
\begin{align*}
L(u_1\,+\,u_2)\,=\,\mathcal {T}_{u_1}(u_2).
\end{align*}
Note that the proof above holds for                                                                     
$h\,=\,x^{2}g\,\in\,C^{4,\,\alpha}(\overline{M})$. Now we have
proved the following theorem,
\begin{thm}\label{LuckFuture}
Let $(M,\,g)$ be an asymptotically hyperbolic manifold of
dimensional $n \geq 4$, with $x$ the smooth defining function, and
the metric $h = x^{2}g\in C^{4,\alpha}(\overline{M})$.                  
For $0\,<\,\nu\,<\,\frac{n-1}{2}$ and $0<\alpha<1$, let                
\begin{align*}
L\,:\,x^{\nu}\,\Lambda^{4,\,\alpha}(M)\,\rightarrow\,x^{\nu}\Lambda^{0,\,\alpha}
\end{align*}
be the linear operator defined in $(\ref{linearoperator})$, which by                                      
Theorem \ref{HLSemiFredholm} is essentially surjective. Assume that
$L$ is surjective. Then there exists a small constant
$\epsilon_0\,>\,0$, depending on the diameter of $M$ with respect to
$h$, $\nu$, $n$ and also $P_1$ and $L$, so that the following holds:

Let $\epsilon$ be any small real number satisfying
$0\,<\,\epsilon\,<\epsilon_0$, and let                                                     
$f\,\in\,\Lambda^{0,\,\alpha}(M)$ satisfy
\begin{align*}
\|Q_g\,-\,f\|_{x^{\nu}\Lambda^{0,\,\alpha}}\,\leq\,\tilde{C}\,\epsilon,
\end{align*}
for some positive constant $\tilde{C}$ depending on the diameter of
$M$ with respect to $h$, $\nu$, $n$, also $P_1$ and $L$.

Then for each $u_1\,\in\,B_{\epsilon}(0)\,\bigcap\,\text{Ker}(L)$,
there exists a unique
$u\,\in\,B_{2\epsilon}(0)\,\subseteq\,x^{\nu}\Lambda^{4,\,\alpha}(M)$,
so that $Q_{\tilde{g}}\,=\,f$, where
$\tilde{g}\,=\,(1\,+\,u)^{\frac{4}{n-4}}g$ for $n\,\geq\,5$, and
$\tilde{g}\,=\,e^{2u}g$ for $n\,=\,4$, with $P_1\,u\,=\,u_1$.
\end{thm}
By the discussion at the end of Section \ref{section_a}, for the
cases in Theorem {\ref{Hyperbolicmetrico}} and Theorem
{\ref{PCPE1}}, $L$ is surjective for
$x^{\nu}\Lambda^{4,\,\alpha}(M),\,\,0\,<\,\nu\,<\,\frac{n\,-\,1}{2}$.           
This completes the proof of \ref{Theorem1_1} of Theorem
\ref{Hyperbolicmetrico} and \ref{PCPE1}.                                        

Since surjectivity is an open property, $L$ is surjective for
$x^{\nu}\Lambda^{4,\,\alpha}(M),\,\,0\,<\,\nu\,<\,\frac{n\,-\,1}{2}$,           
for smooth $g$ that is close enough to these metrics. Theorem
\ref{LuckFuture} holds for metrics in a small neighborhood of these
metrics.

In the following, we will discuss about the boundary
regularity of the solutions. For convenience, we assume that the defining function $x$ and the metric $h\,=\,x^2g$ are smooth up to the boundary. The discussion we use here is standard, see \cite{Mazzeo3}. We will sketch the discussion. Composing the inverse $G$ operator of $L$ on both sides of $(\ref{general_version_nonlinear_equation})$,   
\begin{align}\label{inverse operator_regularity_test}
u\,-\,P_1 u\,=\,G\,L\,u\,=\,G\,\mathcal {T} (u),
\end{align}
with $u_1\,=\,P_1\,u$ the projection of $u$ to the null space of
$L$.


For the regularity of $u$ with respect to the derivative
$\partial_y$, which is the derivative in some $y$ direction, we
introduce the following weighted space with $k\,\leq\,m$:
\begin{align*}
x^{\nu}\Lambda^{m,\,\alpha,\,k}\,=\,\{\,&u\in\,x^{\nu}\Lambda^{m,\,\alpha}(M,\,\sqrt{dx
dy}),\,\,\text{so that
}\,\,(x\partial_x)^j(x\partial_y)^{\beta}\partial_y^{\gamma}u\,\in\,x^{\nu}\Lambda^{0,\,\alpha},\\
&\text{for
}\,\,j\,+\,|\beta|\,+\,|\gamma|\,\leq\,m,\,\,j\geq\,0,\,\,\text{and}\,\,|\gamma|\,\leq\,k.\,\}.
\end{align*}
An easy observation is that for
$u\,\in\,x^{\nu}\Lambda^{m,\,\alpha}$ and $m\,\geq\,1$,
$\partial_yu\,=\,x\partial_y(x^{-1}u)$, so that
\begin{align}\label{regularity_equation_boundary_direction}
\partial_yu\,\in\,x^{\nu-1}\Lambda^{m-1,\,\alpha}.
\end{align}
Also for $u\,\in\,x^{\nu}\Lambda^{m,\,\alpha,\,k}$ and
$1\,\leq\,k\,\leq\,m$,
$\partial_yu\,\in\,x^{\nu}\Lambda^{m-1,\,\alpha,\,k-1}$. In
Proposition 2.9 in \cite{Mazzeo3}, it is proved that the inverse
operator
$G:\,x^{\nu}\Lambda^{m,\,\alpha,\,k}\,\to\,x^{\nu}\Lambda^{m+4,\,\alpha,\,k}$
is bounded for $m\,\geq\,0$ and $0\,\leq\,k\,\leq\,m$; also,
$P_1:\,x^{\nu}\Lambda^{m\,+\,4,\,\alpha,\,k}\,\to\,x^{\nu}\Lambda^{m\,+\,4,\,\alpha,\,k}$
is bounded for $m\,\geq\,0$ and $0\,\leq\,k\,\leq\,m$.

\begin{lem}\label{The induction_lemma}
Let $u\,\in\,x^{\nu}\Lambda^{4,\,\alpha}$ be a solution to
$(\ref{general_version_nonlinear_equation})$ with
$1\,\leq\,\nu\,<\,\frac{n-1}{2}$ and $0<\alpha<1$. Assume that                                           
$(f\,-\,Q_g)\,\in\,x^{\nu}\Lambda^{m,\,\alpha,\,k}$, and
$u_1\,=\,P_1\,u\,\in\,x^{\nu}\Lambda^{m\,+\,4,\,\alpha,\,k}$, for
$0\,\leq\,k\,\leq\,m$. Then we have that
$u\,\in\,x^{\nu}\Lambda^{m+4,\,\alpha,\,k}$.
\end{lem}

\noindent Proof of Lemma \ref{The induction_lemma}. By assumption,
$x$ and the metric $h$ are smooth up to the boundary, so that
$Q_g\,\in\,C^{\infty}(\overline{M})\,\subseteq\,\Lambda^{m,\,\alpha,\,k}$ for any $m\,\geq\,k$,  and then                           
we have $f\,\in\,\Lambda^{m,\,\alpha,\,k}$. For $m\,=\,0$ the claim
holds automatically. Now assume $m\,\geq\,1$. Using $(\ref{inverse
operator_regularity_test})$ and boundedness of $G$ for $k\,=\,0$ we                   
obtain that $u\,\in\,x^{\nu}\Lambda^{1+4,\,\alpha}$. Then we can
substitute the regularity of $u$ into the right hand side of
$(\ref{inverse operator_regularity_test})$, to gain more regularity.
Using this induction argument, we obtain
$u\,\in\,x^{\nu}\Lambda^{m\,+\,4,\,\alpha}\,=\,x^{\nu}\Lambda^{m\,+\,4,\,\alpha,\,0}$.
This proves the lemma for $k\,=\,0$.                                                  

Define the function $F$ on $\mathbb{R}$ as follows,
\begin{align*}
F(u)\,=\,\left\{\begin{array}{l@{\quad\quad}l} e^{4u}\,-\,1\,-\,4u, & n\,=\, 4,\\ \\
(1\,+\,u)^{\frac{n+4}{n-4}}\,-\,1\,-\,\frac{n+4}{n-4}u, &
n\,\geq\,5.\end{array} \right.
\end{align*}
Noticing that for $u\,\in\,x^{\nu}\Lambda^{m,\,\alpha,\,k'}$ with
$k'\,<\,k$, using $(\ref{regularity_equation_boundary_direction})$
and the fact $\nu\,\geq\,1$, we have that
\begin{align*}
u^2f=x\,u\,(x^{-1}u)f\,\in\,x^{\nu}\Lambda^{m,\,\alpha,\,k'+1},
\end{align*}
raising the third index by $1$. This holds for the term $F(u)f$,
since $F$ is smooth on $\mathbb{R}$ and vanishes quadratically at
$0$. Similarly,
\begin{align*}
u(f\,-\,Q_g)\,=\,x\,u\,(x^{-1}(f\,-\,Q_g))\,=\,x\,(x^{-1}u)\,(f\,-\,Q_g)\,\in\,x^{\nu}\Lambda^{m,\,\alpha,\,k'+1}.
\end{align*}
By this fact, combining with the equation $(\ref{inverse
operator_regularity_test})$, and also with boundedness of $G$, an
induction argument as the case $k\,=\,0$ proves the Lemma. \qed                    

Now we assume that $f\,=\,Q_g$. Generally,                                         
$u_1\,=\,P_1\,u\,\in\,x^{\nu}\Lambda^{4,\,\alpha}$ does not have
better regularity. In $(\ref{inverse operator_regularity_test})$,
the terms on the right hand side behave better than $P_1u$, and $u$
behaves like $P_1u$ near the boundary, and $u$ only has the
expansion $(\ref{mainterms})$ with the coefficients which are
distributions of negative order, as discussed in Proposition 3.16 in
\cite{Mazzeo3}. If $1\,\leq\,\nu\,<\,\frac{n-1}{2}$ and
$u_1\,=\,P_1u\,\in\,x^{\nu}\Lambda^{m,\,\alpha,\,k}$ for all                      
$m\,\geq\,k\,\geq\,0$, which as discussed in \cite{Mazzeo1} is
equivalent to say $u_1$ has a smooth expansion $(\ref{sequence0})$,               
then by Lemma \ref{The induction_lemma}, $u$ has a smooth expansion               
as in $(\ref{fomalexpansion})$. Also, for $u_1$ small enough, we
already obtain the existence of $u$ in Poincar$\acute{\text{e}}$
Einstein manifolds. This completes the proof of Theorem
\ref{Hyperbolicmetrico} and Theorem \ref{PCPE1}. \qed
Here we observe that the expansion of $u$ gives us information on                                           
the asymptotic behavior of the curvature. For $n=4$, assume that $g$
and $\tilde{g}$ are asymptotically hyperbolic metrics on $M$, with
the transformation $\tilde{g}\,=\,e^{2u}g$, such that $u$ has the
expansion $u\,\sim\,x^{\frac{3}{2}+
i\frac{\sqrt{15}}{2}}u_{00}(y)\,+\,x^{\frac{3}{2}\,-\,i\frac{\sqrt{15}}{2}}u_{10}(y)\,+\,o(x^{\frac{3}{2}})$.                        
Let $(1\,+\,v)^2\,=\,e^{2u}$. Denote
$\nu_0\,=\,\frac{3}{2}+i\frac{\sqrt{15}}{2}$, and
$\nu_1\,=\,\bar{\nu_0}$. Then,
\begin{align*}
R_{\tilde{g}}\,=\,&(1+v)^{-3}(-6 \Delta_g + R_g)(1+v)\,=\,e^{-3u}(-6
\Delta_g\,+\,R_g)e^u\\
=&-6e^{-u}[-3x\partial_xu\,+\,(x\partial_x)^2u]\,+\,R_g\,-\,2R_gu\,+\,R_g(e^{-2u}-1+2u)\\
&+\,6x^2e^{-3u}(\Delta_y e^u\,+\,\frac{1}{2}\displaystyle\sum_{4\geq
i,j\geq 2}h^{ij}\partial_xe^u).
\end{align*}
Therefore,
\begin{align*}
R_{\tilde{g}}-R_g=&-6e^{-2u}[-3x\partial_xu+(x\partial_x)^2u]-2R_gu+R_g(e^{-2u}-1+2u)+6x^2e^{-3u}(\Delta_ye^u\\              
&+\frac{1}{2}\sum_{4\geq
i,j \geq2}h_{ij}\partial_xh_{ij}\partial_xe^u)\\
=&-6(-3x\partial_xu+(x\partial_x)^2u)+6(1-e^{-2u})(-3x\partial_xu+(x\partial_x)^2u)+24u\\                            
&-2u(12+R_g)
+R_g(e^{-2u}-1+2u)+6x^2e^{-3u}[\Delta_ye^u+\frac{1}{2}\sum_{4\geq
i,j \geq 2}h^{ij}\partial_xh_{ij}\partial_xe^u])\\
=&-6(-3\nu_0x^{\nu_0}u_{00}(y)+\nu_0^2x^{\nu_0}u_{00}(y)\,-\,3\nu_1x^{\nu_1}u_{10}(y)+\nu_1^2x^{\nu_1}u_{10}(y)+O(x^{\frac{3}{2}+1}))\\
&+24(x^{\nu_0}u_{00}(y)+x^{\nu_1}u_{10}(y)+O(x^{\frac{3}{2}+1}))                                                     
+O(x^{\frac{3}{2}+1})\\
=&-6((\nu_0^2-3\nu_0-4)x^{\nu_0}u_{00}(y)+(\nu_1^2-3\nu_1-4)x^{\nu_1}u_{10}(y))+O(x^{\frac{3}{2}+1})\\
=&120u\,+\,o(x^{\frac{3}{2}})                                                                             
\end{align*}
For asymptotically hyperbolic manifolds of higher dimension, with
similar calculation, we obtain the formula                                                                
\begin{align*}
R_{\tilde{g}}\,-\,R_g\,=\,\frac{4(n-1)(n^2+2n-4)}{(n-4)}u\,+\,o(x^{\frac{n-1}{2}}).                       
\end{align*}
\vskip 5pt

\section{Constant $Q$-curvature metrics for perturbed conformal structures}
\vskip10pt

Let $(M,g_0)$ be a Poincar$\acute{\text{e}}$-Einstein manifold,
with
a defining function $x$ and the metric $h_0\,=\,x^2\,g_0$ smooth up to                          
the boundary. Let                                                                            
\begin{align*}
\mathfrak{M}_{\tau}=\{\,h:\,&\text{metrics on}\,\overline{M},                                
\,\text{so that}\,h\in
C^{4,\,\alpha}(\overline{M}),\\
&\text{with}\,\|h-h_0\|_{C^{4,\alpha}(M)} \leq
\tau,\,\text{and}\,|dx|_{h}\big|_{\partial M} = 1\},
\end{align*}
for $\tau\,>\,0$ and $0<\alpha<1$. For                                                     
$h\,\in\,\mathfrak{M}_{\tau}$, let $g\,=\,x^{-2}h$. We want to see
that if $\tau$ is small enough, whether we can find a constant
$Q$-curvature metric $\tilde{g}$ in the conformal class of $g$, with
$Q_{\tilde{g}}\,=\,Q_{g_0}$. We use the same notation $u$, $L_g$ and
so on as above. Note that the
choice of $x$ that $|dx|_{h}\,=\,1$ in the sections before is only                         
to make the notation simpler. Now we only assume that $|dx|_h\,=\,1$
on
$\partial M$, and then there are only some additional small terms in                       
$E(L)$. 
It is easy to check that
\begin{align*}
x^{\alpha}\Lambda^{0,\,\alpha}(M,\,\sqrt{dx\,dy})\,=\,\{\,u\,\in\,C^{\alpha}(\overline{M}),\,u\big|_{\partial
M}\,=\,0\,\}.
\end{align*}
Let $L_{g}$ and $L_{g_0}$ be the linear operators
$(\ref{linearoperator})$ with respect to $g$ and $g_0$. Recall that
$\text{Ric}_g$ and $\text{R}_g$ satisfy $(\ref{Ricciformula})$ and
$(\ref{scalarcurvatureform})$. We know that
\begin{align*}
(|dx|_h^2\,-\,1)\,\in\,x^{\alpha}\Lambda^{0,\,\alpha}(M,\,\sqrt{dx\,dy}),\,\,\text{and}\,\,\|\,|dx|_h^2\,-\,|dx|_{h_0}^2\,\|_{x^{\alpha}\Lambda^{0,\,\alpha}}\,\leq\,C\,\tau,
\end{align*}
for some constant $C$ depending on the defining function and $h_0$.
Also it is easy to see the following inequalities by the formula of
the coefficients
\begin{align*}
&\| (\Delta_g^2\,-\,\Delta_{g_0}^2)u\|_{x^{\alpha}\Lambda^{0,\,\alpha}}\,\leq\,C\,\tau\,\|u\|_{x^{\alpha}\Lambda^{4,\, \alpha}},\\
&\|(R_g\Delta_g\,-\,R_{g_0}\Delta_{g_0})\,u\|_{x^{\alpha}\Lambda^{0,\, \alpha}}\,\leq\,C\,\tau\,\|u\|_{x^{\alpha}\Lambda^{4, \,\alpha}},\\
&\|(\text{Ric}_{ij}(g)\nabla_g^i\nabla_g^j\,-\,\text{Ric}_{ij}(g_0)\nabla_{g_0}^i\nabla_{g_0}^j)\,u\|_{x^{\alpha}\Lambda^{0,\, \alpha}}\,\leq\,C\,\tau\,\|u\|_{x^{\alpha}\Lambda^{4, \,\alpha}},\\
&\|(\nabla_gR_g,\,\nabla_g
u)\,-\,(\nabla_{g_0}R_{g_0},\,\nabla_{g_0}u)\|_{x^{\alpha}\Lambda^{0,\,\alpha}}\,\leq\,C\,\tau\,\|u\|_{x^{\alpha}\Lambda^{4,                    
\,\alpha}},\\
&\|Q_g\,-\,Q_{g_0}\|_{x^{\alpha}\Lambda^{0,\,\alpha}}\,\leq\,C\,\tau,                                      
\end{align*}
with $C$ depending on the defining function $x$ and the metric                                            
$h_0$. We know that $L_{g_0}$ is surjective. Let
\begin{align}
x^{\alpha}\Lambda^{4,\,\alpha}(M,\,\sqrt{dx\,dy})\,=\,\text{Ker}(L_{g_0})\,\oplus\,V_1(g_0),
\end{align}
be the splitting as in Theorem \ref{HLSemiFredholm}. Restricted on                                        
$V_1$ with respect to $g_0$, $L_{g_0}$ satisfies
$(\ref{uniformlybounded})$. Therefore, we can choose $\tau\,>\,0$
small enough so that $\|L_g\,-\,L_{g_0}\|\,\leq\,\frac{1}{2}\,C_1$
with $C_1$ in
$(\ref{uniformlybounded})$. 
Then we have that
\begin{align}
\frac{1}{2}C_1\|u\|_{x^{\alpha}\Lambda^{4,
\,\alpha}}\,\leq\,\|L_g\,u\|_{x^{\alpha}\Lambda^{0,
\,\alpha}}\,\leq\,(C_2\,+\,\frac{1}{2}C_1)\|u\|_{x^{\alpha}\Lambda^{4,\,\alpha}},
\end{align}
for $u\,\in\,V_1(g_0)$. Then
$L_g:\,V_1(g_0)\,\to\,x^{\alpha}\Lambda^{0,
\,\alpha}(M,\,\sqrt{dx\,dy})$ is isomorphic so that
\begin{align}
\|L_g^{-1}\|\,\leq\,\frac{2}{C_1}.
\end{align}
and then $\text{Ker}(L_g)\,\subseteq\,\text{Ker}(L_{g_0})$. We will
only use the splitting of the weighted space with respect to $g_0$.
Now we have a uniform constant $\epsilon
> 0$ for all $h\,\in\,\mathfrak{M}_{\tau}$ and $g\,=\,x^2h$ so that it satisfies the conditions (\ref{condition_ineqn1}) and (\ref{condition_ineqn2}). 
Furthermore, we assume that $\tau >0$ is small enough so that
\begin{align}
\|Q_g\,-\,Q_{g_0}\|_{x^{\alpha}\Lambda^{0,\,\alpha}}\,\leq\,C\,\tau,                                      
\end{align}
and it satisfies corresponding inequalities as $(\ref{closequantity})$ and $(\ref{smallness})$. Therefore, the proof of Theorem        
\ref{LuckFuture} applies. We then obtain the following perturbation
result.
\begin{thm}\label{generating}
Let $(M,g_0)$ be a Poincar$\acute{\text{e}}$-Einstein manifold with
defining function $x$ and the metric $h_0\,=\,x^2\,g_0$ smooth up to                          
the boundary, and let $\mathfrak{M}_{\tau}$ be as above, with $\tau                           
> 0$. There exists $\tau_0\,>\,0$, so that for
$0\,<\,\tau\,<\,\tau_0$, and any
metric $h\,\in\,\mathfrak{M}_{\tau}$, 
there always
exist a family of asymptotically hyperbolic metrics in the conformal
class of $g\,=\,x^{-2}h$ with constant $Q$-curvature $Q_{g_0}$,
which are parametrized by elements in $\text{Ker}(L_{g_0})$.
\end{thm}


 \vskip0.2in

\section {Critical Metrics of Regularized Determinants}
\vskip10pt


Let $M$ be a fourth dimensional asymptotically hyperbolic manifold,
with complete metric $g$ and its smooth defining function $x$, so
that $h\,=\,x^2\,g$ is a smooth metric on $\overline{M}$. Consider the equation    
\begin{align}\label{curvature1_Definition}
U\,=\,U_g\,\equiv\,\gamma_1\,|W|^2\,+\,\gamma_2Q\,-\,\gamma_3\Delta\,R\,=\,C,
\end{align}
where $\gamma_1,\,\gamma_2,\, \gamma_3$ and $C$ are some constants,                       
$W$ is the Weyl tensor, and $Q$, $R$ the $Q$-curvature and the
scalar curvature with respect to $g$. The equation arises as the
Euler-Laglange equation for the regularized determinants,
\begin{align*}
F_A[w]\,=\,\log(\frac{\det A_{\tilde{g}}}{\det A_g}),
\end{align*}
of a conformally covariant operator $A\,=\,A_g$, under the conformal
change of metrics $\tilde{g}\,=\,e^{2w}\,g$, see Chpater 6 in \cite{Andreas Juhl}. More precisely, under the conformal change, 
\begin{align}\label{Euler-Lagrange_equation}
\tilde{U}e^{4\,w}\,=\,&U\,+\,(\frac{1}{2}\gamma_2\,+\,6\,\gamma_3)\Delta^2w\,+\,6\gamma_3\Delta\,|\nabla\,w|^2\,-\,12\,\gamma_3\,\nabla^i[(\Delta\,w\,+\,|\nabla\,w|^2)\nabla_iw]\\
                       &+\,\gamma_2R_{ij}\nabla_i\nabla_jw\,+\,(2\gamma_3\,-\,\frac{1}{3}\gamma_2)R\Delta\,w\,+\,(2\,\gamma_3\,+\,\frac{1}{6}\gamma_2)(\nabla\,R,\,\nabla\,w),    
\end{align}
with $\tilde{U}\,=\,U_{\tilde{g}}$. 
Define $\alpha\,=\,\frac{\gamma_2}{12\gamma_3}$. The following are
some examples that we are interested in. 
\begin{exm}
For the conformal Laplacian, $A=L$, we have that
$(\gamma_1,\gamma_2,\gamma_3)=(1,\,-4,\,-\frac{2}{3})$, and                                     
$\alpha=\frac{1}{2}$.
\end{exm}
\begin{exm}
For the spin Laplacian, $A\,=\,D^2$, we  have that
$(\gamma_1,\,\gamma_2,\,\gamma_3)\,=\,(7,\,-88,\,-\,\frac{14}{3})$,
and $\alpha\,=\,\frac{11}{7}$.
\end{exm}
\begin{exm}
For the Paneitz operator, $A\,=\,P$, we have that
$(\gamma_1,\,\gamma_2,\,\gamma_3)\,=\,(-\frac{1}{4},\,-14,\,\frac{8}{3})$,
and $\alpha\,=\,\frac{-7}{16}$.
\end{exm}

For convenience, dividing both sides of the function by
$6\,\gamma_3$, we have the following equation,
\begin{align}\label{Equationof_U-curvature}
\frac{\tilde{U}}{6\gamma_3}e^{4w}\,=\,&(1\,+\,\alpha)\,\Delta^2w\,+\,\Delta|\nabla
w|^2\,-\,2\,\nabla^i[(\Delta\,w\,+\,|\nabla\,w|^2)\nabla_iw]\,+\,2\alpha\,R_{ij}\nabla^i\nabla^j
w\\
& +\,(\frac{1}{3}\,-\,\frac{2}{3}\alpha)R\Delta
w\,+\,(\frac{1}{3}\,+\,\frac{1}{3}\alpha)(\nabla R, \nabla
w)\,+\,\frac{U}{6\gamma_3}.                                                       
\end{align}
We should note that                                                                               
\begin{align*}
\Delta |\nabla w|^2\,-\,2\nabla^i(\Delta w\,\nabla_i
w)\,=\,&2(\Delta_g\nabla w,\,\nabla w)\,+\,2(\nabla^2 w,\,\nabla^2                            
w)\,-\,2\,\nabla^i(\Delta w\,\nabla_i w)\\
=\,&2\nabla^i w(g^{pq}\nabla_p\nabla_q\nabla_i
w-g^{pq}\nabla_i\nabla_p\nabla_q                                                              
w)+2(|\nabla^2w|_g^2-(\Delta w)^2)\\                                                          
=\,&2\nabla^i w
g^{pq}R_{piq}^s\nabla_s
w\,+\,2(|\nabla^2w|_g^2\,-\,(\Delta w)^2)\\
=\,&2\,Ric(\nabla w,\,\nabla w)\,+\,2(|\nabla^2 w|_g^2\,-\,(\Delta
w)^2).
\end{align*}
Moreover,
\begin{align*}
\nabla^i(|\nabla w|^2 \nabla_i
w)\,=\,2\,\nabla^i\nabla^jw\,\nabla_jw\nabla_iw\,+\,|\nabla
w|^2\,\Delta w,
\end{align*}
therefore, the equation can be written in the following way,
\begin{align}\label{nonlinearequations}
\frac{\tilde{U}}{6\gamma_3}e^{4w}\,=\,&(1\,+\,\alpha)\,\Delta^2w\,+\,2\,Ric(\nabla
w, \,\nabla w)\,+2(|\nabla^2 w|_g^2\,-\,(\Delta
w)^2)\\
&-\,4\,\nabla_i\nabla_jw\,\nabla_jw\nabla_iw                                                
-\,2\,|\nabla w|^2\,\Delta w\,+\,2\alpha\,R_{ij}\nabla^i\nabla^j                            
w\\
&+\,(\frac{1}{3}\,-\,\frac{2}{3}\alpha)R\Delta
w\,+\,(\frac{1}{3}\,+\,\frac{1}{3}\alpha)(\nabla R, \nabla
w)\,+\,\frac{U}{6\gamma_3}.                                                                 
\end{align}
 We should point out that for $\alpha\,=\,-1$ and $\gamma_1\,=\,0$, the equation reduces to a second order differential equation, and in this case the $U$-curvature relates to the $\sigma_2$-curvature with respect to the Schouten tensor  
 $A(g)$,
\begin{align*}
\frac{1}{12\gamma_3}U(g)\,&=\,\frac{\gamma_1}{12\gamma_3}|W|_g^2\,+\,\frac{\gamma_2}{12\gamma_3}Q_g\,-\,\frac{\Delta R_g}{12}\\
&=\,-(\frac{-1}{4}|\text{Ric}_g|^2\,+\,\frac{1}{12}R_g^2\,-\,\frac{1}{12}\Delta_gR_g)\,-\,\frac{\Delta R_g}{12}\\
&=\,-\,(\frac{-1}{4}|\text{Ric}_g|^2\,+\,\frac{1}{12}R_g^2)\,=\,-\,2\sigma_2(g).                                        
\end{align*}
We have the equation
\begin{align*}
4\sigma_2(\tilde{g})\,=\,&-\,2\text{Ric}(\nabla_g w,\,\nabla_g w)\,-\,2(|\nabla^2w|_g^2\,-\,(\Delta w)^2)\,+\,4\,\nabla_{ij}w\nabla_iw\nabla_jw\\
&+\,2|\nabla w|^2\Delta w\,+\,2
\text{Ric}_{ij}\nabla^i\nabla^jw\,-\,R_g \Delta w\,+\,4\sigma_2(g).                                                  
\end{align*}
A prescribed constant $\sigma_2$-curvature asymptotically hyperbolic metric problem is                                                                       
discussed in \cite{MP}. From now on, we assume that
$\alpha\,\neq\,-1$.

The linearization of $(\ref{nonlinearequations})$ is given by                                                       
\begin{align*}
L\,w=(1+\alpha)\Delta^2w+2\alpha R_{ij}\nabla^i\nabla^j                                                     
w+(\frac{1}{3}-\frac{2}{3}\alpha)R\Delta
w+(\frac{1}{3}+\frac{1}{3}\alpha)(\nabla R, \nabla w) - \frac{2 U}{3                                        
\gamma_3}w=0.
\end{align*}
As $x\,\to\,0$,                                                                                                     
\begin{align*}
R_{ijkl}(g)\,=\,&x^{-2}[R_{ijkl}(h)\,-\,h_{ik}(x^{-1}\nabla_j^h\nabla_l\,x\,+\,\frac{1}{2}x^{-2}h_{jl})\,-\,h_{jl}(-x^{-1}\nabla_i^h\nabla_k\,x\,+\,\frac{1}{2}x^{-2}h_{ik})\\
&+\,h_{il}(- x^{-1}\nabla_j^h\nabla_k
x\,+\,\frac{1}{2}x^{-2}h_{jk})\,+\,h_{jk}(-x^{-1}\nabla_i^h\nabla_l
x\,+\,\frac{1}{2}x^{-2}h_{il})]\\
=\,&x^{-4}[-\frac{1}{2}h_{ik}h_{jl}\,-\,\frac{1}{2}h_{jl}h_{ik}\,+\,h_{il}h_{jk}\,+\,\frac{1}{2}h_{jk}h_{il}\,+\,O(x)]\\
=\,&x^{-4}[-\,h_{ik}h_{jl}\,+\,h_{il}h_{jk}\,+\,O(x)],
\end{align*}
while
\begin{align*}
A(g)=\frac{1}{4-2}(Ric(g)-\frac{1}{2(4-1)}R(g)g)=\frac{1}{2}(-3+2+O(x))g=(-\frac{1}{2}+O(x))g,            
\end{align*}
so that
\begin{align*}
W_{ijkl}(g)\,=\,&R_{ijkl}(g)\,-\,g_{ik}A_{jl}(g)\,+\,g_{il}A_{jk}(g)\,+\,g_{jk}A_{il}(g)\,-\,g_{jl}A_{ik}(g)\\
=\,&x^{-4}(-\,h_{ik}h_{jl}\,+\,h_{il}h_{jk}\,+\,O(x))\,+\,x^{-4}[-
h_{ik}(-\frac{1}{2}h_{jl}\,+\,O(x))\\
&+\,h_{il}(-\frac{1}{2}h_{jk}\,+\,O(x))\,+\,h_{jk}(-\frac{1}{2}h_{il}\,+\,O(x))\,-\,h_{jl}(-\frac{1}{2}h_{ik}\,+\,O(x))]\\
=\,&x^{-4}O(x),
\end{align*}
and moreover, using the fact $\Delta_hR\,=\,O(x)$, and
$Q(g)\,=\,3\,+\,O(x)$, we have that $U(g)\,=\,3\gamma_2\,+\,O(x)$.
Then we obtain the main terms of $L\,w$ as follows,                         
\begin{align*}
L\,w=&\,(1+\alpha)\Delta_g^2w+(\frac{1}{3}-\frac{2}{3}\alpha)R_g\Delta_gw+2\alpha\,\text{Ric}_{ij}^g\nabla_g^i\nabla_g^j                       
w\\
&+\frac{1}{3}(1+\alpha)(\nabla_gR_g,\,\nabla_gw)-\frac{2
U}{3\gamma_3}w\\                                                                                                                               
=\,&(1+\alpha)\Delta_g^2w+(\frac{1}{3}-\frac{2}{3}\alpha)(-12+O(x))\Delta_gw+2\alpha(-3\Delta_gw+O(x)p(x,y,x\partial_x, x\partial_y)w)\\                
&+\frac{1}{3}(1\,+\,\alpha)(-(2\times
4\,-\,2)x^2H(h|_{S_x})\partial_xw\,+\,O(x^3)|\nabla_yw|)-\,(8\,\times
3\,\alpha\,+\,O(x))w\\
=\,&(1\,+\,\alpha)\Delta_g^2w-12(\frac{1}{3}-\frac{2}{3}\alpha)\Delta_gw-6\alpha\Delta_g
w - 24\,\alpha\,w+O(x)p(x,\,y,\,x\partial_x,\,x\partial_y)w
\\
=\,&(1\,+\,\alpha)\Delta_g^2w-(4\,-\,2\alpha)\Delta_gw\,-\,24\,\alpha\,w\,+\,O(x)p(x,\,y,\,x\partial_x,\,x\partial_y)w\\              
=\,&((1\,+\,\alpha)\Delta_g\,+\,6\,\alpha)(\Delta_g\,\,-\,4)w\,+\,O(x)p(x,\,y,\,x\partial_x,\,x\partial_y)w.                              
\end{align*}
Correspondingly,                                                                                          
\begin{align*}
N(L)w&=(1+\alpha)((s\partial_s)^2+s^2\Delta_v-3s\partial_s)^2w+(2\alpha-4)((s\partial_s)^2+s^2\Delta_v-3s\partial_s)w-24\alpha w,
\end{align*}
\begin{align*}
L_0(t,\,\hat{\eta})w&=(1+\alpha)((t\partial_t)^2+t^2-3t\partial_t)^2w+(2\alpha-4)((t\partial_t)^2+t^2-3t\partial_t)w-24\alpha w\\
  &=((1+\alpha)((t\partial_t)^2+t^2-3t\partial_t)+6\alpha
       )(((t\partial_t)^2+t^2-3t\partial_t)-4)w=L_3\circ\,L_1 w,\\
       \\
I(L)w&=\,(1\,+\,\alpha)((s\partial_s)^2\,-\,3\,s\partial_s)^2w\,+\,(2\,\alpha\,-\,4)((s\partial_s)^2\,-\,3s\partial_s)w\,-\,24\,\alpha\,w\\           
      &=\,((1\,+\,\alpha)((s\partial_s)^2\,-\,3\,s\partial_s\,+\,6\alpha)((s\partial_s)^2\,-\,3\,s\partial_s\,-\,4)w.                                   
\end{align*}
Therefore, the indicial roots of $L$ is as follows,
\begin{enumerate}[label=\roman{*})]                                                                    
\item For $\alpha\,=\,\frac{1}{2}$,
$spec_b(L)\,=\,\{4,\,-1,\,1,\,2\}$.

\item For $\alpha\,=\,\frac{11}{7}$,
$spec_b(L)\,=\,\{4,\,-1,\,\frac{3}{2}\,+\,i\frac{\sqrt{51}}{6},\,\frac{3}{2}\,-\,i\,\frac{\sqrt{51}}{6}\}$.

\item For $\alpha\,=\,\frac{-7}{16}$,
$spec_b(L)\,=\,\{4,\,-1,\,\frac{3}{2}\,+\,\frac{\sqrt{249}}{6},\,\frac{3}{2}\,-\,\frac{\sqrt{249}}{6}\}$.
\end{enumerate}
The solution of $L_1 w =0$ is exactly the same as discussed in
Section $2$. We solve $L_3\,w\,=\,0$ by transferring it into the               
Bessel type equations discussed as above. Let
$u(t)\,=\,t^{\beta}\tilde{w}(t)$, then                            
\begin{align*}
0\,=\,&t^{\beta}((t\partial_t)^2\tilde{w}\,+\,(2\beta\,-\,3)t\partial_t\tilde{w}\,+\,(\beta^2\,-\,3\beta\,+\frac{6\alpha}{1\,+\,\alpha}\,-t^2)\tilde{w}).         
\end{align*}
Let $2\beta\,-\,3\,=\,0$, and then $\beta\,=\,\frac{3}{2}$.
Consequently,
\begin{align*}
[(t\partial_t)^2\,-\,(t^2\,+\,\frac{9}{4}\,-\,\frac{6\alpha}{1+\alpha})]\tilde{w}\,=\,0.            
\end{align*}
Let $\tilde{\alpha}^2\,=\,\frac{9}{4}\,-\,\frac{6
\alpha}{1+\alpha}$, then the solution is
\begin{align}\label{2}
w\,=\,t^{\frac{3}{2}}(C_1I_{\tilde{\alpha}}(t)\,+\,C_2K_{\tilde{\alpha}}(t)).                             
\end{align}
Here $\tilde{\alpha}^2$ is
$\frac{1}{4},\,\frac{-17}{12},\,\frac{83}{12}$, corresponding to the                                         
above three cases, with $\text{Re}(\tilde{\alpha})\,\geq\,0$. For
the case $\tilde{\alpha}^2\,=\,-\,\frac{17}{12}$, since                                                   
$\tilde{\alpha}^2$ is negative, $L_3$ behaves the same as $L_2$ in
Section $2$, and it follows that Theorem \ref{SLSemiFredholm} and                                         
Theorem \ref{HLSemiFredholm} with $n\,=\,4$ hold for the linear
operator $L$, using the same argument as in Section $2$.

By the expansion of the series form of the Bessel functions, as in
[\cite{N.N.L}, P. 108], we have
\begin{align*}
t^{\frac{3}{2}}I_{\tilde{\alpha}}(t|\eta|)\sim\,t^{\frac{3}{2}+\tilde{\alpha}}|\eta|^{\tilde{\alpha}}/(2^{\tilde{\alpha}}\Gamma(1\,+\,\tilde{\alpha})),
\end{align*}
and
\begin{align*}
t^{\frac{3}{2}}I_{-\,\tilde{\alpha}}(t|\eta|)\sim\,t^{\frac{3}{2}-\tilde{\alpha}}|\eta|^{-\tilde{\alpha}}/(2^{-\tilde{\alpha}}\Gamma(1\,-\,\tilde{\alpha})),
\end{align*}
 near $t\,=\,0$. Here we should note that the series expansion applies for all $\tilde{\alpha}\,\in\,\mathbb{C}$. Now it is easy to see that the linear 
 combination
 \begin{align*}
 x^{\frac{3}{2}}(C_1\,x^{\tilde{\alpha}}\,+\,C_2\,x^{-\tilde{\alpha}})
 \end{align*}
 can never vanish to infinite order at $t\,=\,0$ if either $C_1 \neq 0$ or $C_2 \neq 0$. Also,
\begin{align*}
t^{\frac{3}{2}}K_{\tilde{\alpha}}(t|\eta|)\sim\,t^{\frac{3}{2}}\,\frac{\pi}{2}\frac{I_{\tilde{\alpha}}(t|\eta|)\,-\,I_{-\,\tilde{\alpha}}(t|\eta|)}{\sin(\tilde{\alpha}\pi)}\,\sim\,O((t|\eta|)^{\frac{3}{2}-\tilde{\alpha}}),
\end{align*}
near $t\,=\,0$, with $\tilde{\alpha}\,>\,0$ and $\tilde{\alpha}\,\neq\,1,\,2,\,3,...$                             

Using the integral form, we have
\begin{align*}
t^{\frac{3}{2}}I_{\tilde{\alpha}}(t|\eta|)\,\,\text{grows
exponentially
},\,\,t^{\frac{3}{2}}K_{\tilde{\alpha}}(t|\eta|)\,\,\text{decays                                                  
exponentially}
\end{align*}
near $t\,=\,+\,\infty$. 
Therefore, $t^{\frac{3}{2}}I_{\tilde{\alpha}}(t|\eta|)$ does not                                                  
belong to $t^{\delta}L^2(\mathbb{R}^{+})$ for any $\delta > 0$,
while
\begin{align*}
t^{\frac{3}{2}}K_{\tilde{\alpha}}(t|\eta|)\,\in\,t^{\delta}L^2(\mathbb{R}_{+}),
\end{align*}
only for
$\delta\,<\,\frac{3}{2}+\frac{1}{2}\,-\,\tilde{\alpha}\,=\,2\,-\,\tilde{\alpha}$.
That is, $L_3$ is injective in $x^{\delta}L^2$ for                                                               
$\delta\,>\,2\,-\,\tilde{\alpha}$.

Summarizing the above discussion, let us compute $\bar{\delta}$ and $\underline{\delta}$ for the linearized operator $L$.                   
\begin{align*}
&\bar{\delta}\,=\,\text{inf}\,\{\delta:\,L_1\,\,\text{and}\,\,L_3\,\,\text{are
injective
in}\,\,t^{\delta}L^2\}\,=\,\text{sup}\,\{-1\,+\,\frac{1}{2},\,2\,-\,\tilde{\alpha}\},\,\,\text{and dually,}\\
&\underline{\delta}\,=\,\text{inf}\,\{(\frac{3}{2}+\frac{1}{2})\times
2\,-\,(-1\,+\,\frac{1}{2}),\,(\frac{3}{2}+\frac{1}{2})\times                                                                           
2\,-\,(2\,-\,\tilde{\alpha})\}\,=\,\text{inf}\,\{\frac{9}{2},\,2+\tilde{\alpha}\}.
\end{align*}
For the case $\alpha\,=\,\frac{1}{2}$, $\bar{\delta}\,=\,\frac{3}{2}$, and                                                             
$\underline{\delta}\,=\,\frac{5}{2}$( surjectivity). For the case
$\alpha\,=\,-\frac{7}{16}$,
$\bar{\delta}\,=\,-1\,+\,\frac{1}{2}\,=\,-\,\frac{1}{2}$, and
$\underline{\delta}\,=\,\frac{9}{2}$. Then we can use Theorem
\ref{SLSemiFredholm} and Theorem \ref{HLSemiFredholm}, to obtain the
semi-Fredholm property for these linear operators.                                                                            

For the Poincar$\acute{e}$ Einstein manifold $(M,\,g)$, we have that
the $U$ curvatures defined above are all constants on $M$. 
We want to see the solutions of the nonlinear problem. Now
$L\,w\,=\,((1\,+\,\alpha)\Delta_g\,+\,6\,\alpha)(\Delta_g\,\,-\,4)w$.              
Define the operator $\mathcal
{T}\,:\,x^{\nu}\Lambda^{4,\,\alpha}(M)\,\rightarrow\,x^{\nu}\Lambda^{0,\,\alpha}(M)$
as follows,
\begin{align*}
\mathcal {T}(w)\,=\,&(\frac{\tilde{U}}{6\gamma_3}e^{4w}\,-\,\frac{U}{6\gamma_3}\,-\,\frac{2}{3\gamma_3}U\,w)\,-\,2\text{Ric}(\nabla w, \,\nabla w)\\
&-\,2(|\nabla^2w|_g^2\,-\,(\Delta
w)^2)\,+\,4\nabla_j\nabla_iw\nabla^jw\nabla^iw\,+\,2|\nabla                            
w|^2\Delta w.
\end{align*}
We rewrite it in the form
\begin{align*}
\mathcal
{T}(w)\,=\,&\frac{\tilde{U}}{6\gamma_3}(e^{4w}\,-\,1\,-\,4\,w)\,+\,(\tilde{U}\,-\,U)(\frac{1}{6\gamma_3}\,+\,\frac{2}{3\gamma_3}\,w)\,-\,2\text{Ric}(\nabla w, \,\nabla w)\\   
&-\,2(|\nabla^2w|_g^2\,-\,(\Delta
w)^2)\,+\,4\nabla_j\nabla_iw\nabla^jw\nabla^iw\,+\,2|\nabla                            
w|^2\Delta w.
\end{align*}
In this formula, comparing with the nonlinear term defined for
$Q$-curvature equation, a few square terms of $w$ and its                     
derivatives of order up to $2$ are involved, which are small terms
in the argument of the perturbation problem. Now, the nonlinear                                     
equation becomes
\begin{align*}
L_g\,w\,=\,\mathcal {T}(w).                                                                                               
\end{align*}
To solve this, the argument follows exactly the way in Section $3$                                 
and Section $4$. We only need to choose the right weighted
H$\ddot{\text{o}}$lder spaces.
Note that the index of the weight for the H$\ddot{\text{o}}$lder                                                              
space is $\frac{1}{2}$ less than the index of the weight of the
corresponding Sobolev spaces.

\subsection{Summary.}

Perturbation results for the curvatures defined in $(\ref{curvature1_Definition})$ can be proved along the same lines as the $Q$-curvature.              
For instance, assume $(M,\,g)$ is a Poincar$\acute{\text{e}}$-Einstein manifold. For the case $\alpha\,=\,-\frac{7}{16}$, by maximal                         
principle, $((1\,+\,\alpha)\Delta_g\,+\,6\,\alpha)$ and                                                                    
$(\Delta_g\,\,-\,4)$ are both injective on $L^2(M,\,g)$. Then
similar to the discussion for the $Q$-curvature equation, there are                                                        
infinitely many solutions
$u\,\in\,x^{\nu}\Lambda^{4,\,\beta}(M,\,\sqrt{dxdy})$ for
$0<\beta<1$ to this equation parametrized by the projection $P_1u$                               
to the kernel of the
linearized operator $L$, for $\nu\,\in\,(0,\,\frac{3}{2})$. Moreover, if                         
$\tilde{U}\,=\,U$, then $w$ has the weak expansion
$w(x,\,y)\,\sim\,w_{00}(y)x^4\,+\,o(x^4)$, and also $w$ has a smooth expansion if $1\,\leq\,\nu\,<\,\frac{3}{2}$ and $P_1w$ has a smooth expansion.      
For the case $\alpha\,=\,\frac{11}{7}$, it is the same as the $Q$                            
curvature problem, and the only difference is that here we use
$i\,\sqrt{51}$ in the indicial roots and in the formula of expansion                             
to replace $i\,\sqrt{15}$. For the case
$\alpha\,=\,\frac{1}{2}$, $((1\,+\,\alpha)\Delta_g\,+\,6\,\alpha)$ is                            
essentially injective on
$x^{\nu}\Lambda^{4,\,\beta}(M,\,\sqrt{dxdy})$ for $\nu\, >\,1$ and $\nu\neq 2$,                                  
while it is essentially surjective  on
$x^{\nu}\Lambda^{4,\,\beta}(M,\,\sqrt{dxdy})$ for $\nu\,<\,2$, also
$\nu\neq 1$ and $0<\beta<1$. Since $(\frac{3}{2}\Delta_g\,+\,3)$ may have finite                                 
dimensional kernel, we do not have perturbation result for $\nu$ in
this
interval. But note that, using the same argument as in Lemma                                     
\ref{injectivitylem} in weighted H$\ddot{\text{o}}$lder spaces, for
$\nu\,>\,2$, the operator
\begin{align*}
(\frac{3}{2}\Delta\,+\,3):\,x^{\nu}\Lambda^{2\,+\,m,\,\beta}\,\to\,x^{\nu}\Lambda^{m,\,\beta},
\end{align*}
is injective, for $0<\beta<1$ and $m\,\geq\,0$. Then dually the                                  
operator $(\frac{3}{2}\Delta\,+\,3)$ is surjective for
$\nu\,\in\,(0,\,1)$. Also we know that the operator
$(\Delta_g\,-\,4)$ is surjective in the weighted
H$\ddot{\text{o}}$lder space with $0\,<\,\nu\,<\,\frac{3}{2}$, then
the linearized operator
\begin{align*}
L:\,x^{\nu}\Lambda^{4\,+\,m,\,\beta}\,\to\,x^{\nu}\Lambda^{m,\,\beta},
\end{align*}
with $m\,\geq\,0$ is surjective for $0\,<\,\nu\,<\,1$ and
$0<\beta<1$. Therefore, for the case $\alpha\,=\,\frac{1}{2}$, the                               
existence result as in \ref{Theorem2_1} in Theorem \ref{PCPE1} holds
for $0\,<\nu\,<\,1$.
For the boundary expansion when $\tilde{U}\,=\,U$, since all the                                                        
indicial roots are integers in this case, there may be $\log(x)$                                  
terms in the expansion as $(\ref{generalformalexpansion})$. Also,
since $\nu\,<\,1$, the smooth expansion result does not hold.



\vskip 0.5in



\section*{Acknowledgements}
It is my pleasure to express my sincere gratitude to my advisor
Professor Matthew Gursky for introducing this question, many helpful
discussions and great patience in the course of this work. I am also
indebted to Professor Rafe Mazzeo for inviting me to Stanford
University to explain his theory of wedge operators, his many
comments on preliminary versions of the manuscript and his continued
interest in this work. Thanks are also due to Professor C. Robin
Graham for useful discussions and introducing me his paper
\cite{Graham and Zworski}. I would also like to thank Professor
Huicheng Yin in Nanjing University, who introduced the PDE area to
me five years ago and taught me the theory of elliptic operators,
for his continuous interest in this problem and constant support. At
last, but not the least, I would like to thank my friends Ye Li and
Yueh-Ju Lin for their encouragement and constant support.






\begin{thebibliography}{s2}

\bibitem{ACF} L.Andersson, Piotr T. Chru$\acute{\text{s}}$ciel, H. Friedrich, {\em On the regularity of solutions to the Yamabe equation and the existence of smooth hyperboloidal initial data for Einstein's field equations}, Comm. Math. Phys. {\bf{149}}
(1992), no. 3, 587 - 612.

\bibitem{S. Brendle} S. Brendle,  {\em Global existence and convergence for a higher order flow in conformal geometry}, Ann. of Math. {\bf{158}} (2003), 323 -
343.  

\bibitem{Chang-Yang} A. Chang, P. Yang, {\em Extremal metrics of zeta functional determinants on $4$-
manifolds}, Ann. of Math. {\bf{142}} (1995), 171 - 212.

\bibitem{Chen and Xu} X. Chen, X. Xu, {\em $Q$-curvature flow on the standard sphere of even dimension}, Journal of Functional Analysis {\bf{261}} (2011), 934 - 980.

\bibitem{DHL} Z. Djadli, E. Hebey, M. Ledoux, {\em Paneitz-type operators and applications}, Duke Math. J. {\bf {104}} (2000), no. 1,
129¨C169.

\bibitem{DjadliaMalchiodi} Z. Djadli, A. Malchiodi,  {\em Existence of conformal metrics with constant $Q$-curvature}, Ann. of Math. {\bf{168}} (2008), 813 -
858.   

\bibitem{Graham}  R. Graham, {\em Volume and Area Renormalizations for Conformally Compact Einstein
Metrics}, arXiv:math/9909042v1 [math.DG] 8 Sep 1999, preprint.

\bibitem{Graham and Zworski} R. Graham, M. Zworski, {\em Scattering matrix in conformal geometry}, Invent. Math. {\bf{152}} (2003) no. 1, 89 - 118.

\bibitem{GRA}  H. Grunau, M. Ould Ahmedou, M. Reichel, {\em The Paneitz equation in hyperbolic space}, Ann. Inst. H. Poincar$\acute{\text{e}}$ Anal. Non Lin$\acute{\text{e}}$aire, {\bf{25}} (2007), 847 - 864.

\bibitem{Matt.Gursky}  M. Gursky, {\em Weyl functional, de Rham cohomology, and Kahler-Einstein metrics}, Ann. of Math. {\bf{148}}  (1998), 315 - 337.

\bibitem{M.Gursky}  M. Gursky, {\em The principal eigenvalue of a conformally invariant differential operator, with
an application to semilinear elliptic PDE}, Comm. Math. Phys.
{\bf{207}} (1999), 131 - 143.

\bibitem{Andreas Juhl}  A. Juhl, {\em Families of conformally covariant differential operators, $Q$-curvature and holography}, Progress in Mathematics, {\bf{275}}. Birkh$\ddot{\text{a}}$user Verlag, Basel, 2009.

\bibitem{N.N.L}  N. N. Lebedev, {\em Special functions and their applications}, Dover, New York, 1972.

\bibitem{Lin} C. Lin,  {\em A classification of solutions of a conformally invariant fourth order equation in $\mathbb{R}^n$},Comment. Math. Helv. {\bf{73}} (1998), 206 - 231.

\bibitem{Mazzeo1}  R. Mazzeo,  {\em Elliptic Theory of Differential Edge Operators I}, Comm. Partial Differential Equations {\bf{16}}(10) (1991), 1615 - 1664.

\bibitem{Mazzeo2}  R. Mazzeo,  {\em The Hodge Cohomology of A Conformally Compact Metric}, J. Diff. Geo. {\bf{28}} (1988), 309 - 339.

\bibitem{Mazzeo3}  R. Mazzeo,  {\em Regularity for the singular Yamabe problem}, Indiana Univ. Math. J. {\bf{40}} (1991), no. 4, 1277 - 1299.

\bibitem{Mazzeo}  R. Mazzeo,  {\em Unique Continuation at Infinity and Embedded Eigenvalues for Asymptotically Hyperbolic Manifolds}, American Journal of
Mathematics {\bf{113}} (1991), 25 - 45.

\bibitem{MP}  R. Mazzeo, F. Pacard, {\em Poincar$\acute{\text{e}}$-Einstein metrics and the Schouten tensor}, Pacific J. Math. {\bf{212}} no. 1 (2003), 169 - 185.

\bibitem{Ndiay1}  C. Ndiaye,  {\em Constant $Q$-curvature metrics in arbitrary dimension}, Journal of Functional Analysis {\bf{251}} (2007), 1 - 58.

\bibitem{Ndiay}  C. Ndiaye,  {\em Conformal metrics with constant $Q$-curvature for manifolds with boundary}, communications in analysis and geometry {\bf{16}} no. 5 (2008), 1049 - 1124.


\bibitem{GW}  G. N. Watson,  {\em A Treatise on the Theory of Bessel Functions}, Cambridge University Press, Cambridge, 1922.

\bibitem{Wei and Ye} J. Wei, D. Ye,  {\em Nonradial solutions for a conformally invariant fourth order equation in $\mathbb{R}^4$}, Calc. Var. {\bf{32}} no. 2 (2008), 373 - 386.

\bibitem{Xu-Yang}  X. Xu, P. Yang, {\em Positivity of Paneitz Operators}, Discrete and Continuous Dynamical Systems {\bf{7}} no. 2 (2001), 329 - 342.
\end{thebibliography}
\end{document}